\documentclass[12pt,reqno]{amsart}
\usepackage{a4wide}

\numberwithin{equation}{section}

\newtheorem{theorem}{Theorem}[section]
\newtheorem{proposition}[theorem]{Proposition}

\newtheorem{lemma}[theorem]{Lemma}

\theoremstyle{definition}

\theoremstyle{remark}
\newtheorem{remark}[theorem]{Remark}

\newcommand{\ep}{\varepsilon}
\newcommand{\Om}{\Omega}

\newcommand{\de}{\delta}

\begin{document}

\title[  Regularization of point vortices
 ]
{
   Regularization of point vortices  for the Euler equation in dimension
   two, part II
 }
 \author{Daomin Cao
}

\address{Institute of Applied Mathematics, Chinese Academy of Science, Beijing 100190, P.R. China}

\email{dmcao@amt.ac.cn}

 \author{Zhongyuan Liu
}

\address{Institute of Applied Mathematics, Chinese Academy of Science, Beijing 100190, P.R. China}

\email{liuzy@amss.ac.cn}

\author{Juncheng Wei
}

\address{Department of Mathematics, The Chinese University of Hong Kong, Shatin, N.T., Hong Kong }

\email{wei@math.cuhk.edu.hk}

\begin{abstract}
In this paper, we continue to construct  stationary classical
solutions of the incompressible Euler equation approximating
singular stationary solutions of this equation.
 This procedure now is carried out by constructing solutions to the
 following elliptic problem
\[
 \begin{cases}
-\ep^2 \Delta
u=\sum_{i=1}^m\chi_{\Om_i^+}\big(u-q-\frac{\kappa_i^+}{2\pi}\ln\frac{1}{\ep}\big)_+^p
-\sum_{j=1}^n\chi_{\Om_j^-}\big(q-\frac{\kappa_j^-}{2\pi}\ln\frac{1}{\ep}-u\big)_+^p, \quad & x\in\Omega, \\
u=0, \quad & x\in\partial\Omega,
\end{cases}
\]
where $p>1$,  $\Omega\subset\mathbb{R}^2$ is a simply connected
bounded domain, $\Om_i^+$ and $\Om_j^-$ are mutually disjoint
subdomains of $\Omega$, $q$ is a harmonic function.

We showed that if $\Omega$ is a simply-connected smooth domain, then
for any given $C^1$-stable critical point of Kirchhoff-Routh
function $\mathcal{W}(x_1^+,\cdots,x_m^+,x_1^-,\cdots,x_n^-)$ with
$\kappa^+_i>0\,(i=1,\cdots,m)$ and $\kappa^-_j>0\,(j=1,\cdots,n)$,
then there is a
 stationary classical solution approximating stationary $m+n$ points vortex solution of incompressible Euler
equations with total vorticity $\sum_{i=1}^m\kappa_i-\sum_{j=1}^n\kappa_j^-$. \\[12pt]
\emph{AMS 2000 Subject Classifications: Primary $35\mathrm{J}60$;\
Secondary   $35\mathrm{JB}05$;  $35\mathrm{J}40$
      }
\end{abstract}

\maketitle

\section{Introduction and main results}

The incompressible Euler equations
\begin{equation}\label{1.2}
\begin{cases}
 \mathbf{v}_t+(\mathbf{v}\cdot\nabla)\mathbf{v}=-\nabla P,\\
\nabla\cdot\mathbf{v}=0,
\end{cases}
\end{equation}
describe the evolution of the velocity $\mathbf{v}$ and the pressure
$P$ in an incompressible flow. In $\mathbb{R}^2$, the vorticity of
the flow is defined by $\omega=\nabla\times\mathbf{v}:=\partial_1
v_2-\partial_2 v_1$, which satisfies the equation
\[
\omega_t+\mathbf{v}\cdot\nabla\omega=0.
\]

Suppose that $\omega$ is known, then the velocity $\mathbf{v}$ can
be recovered by Biot-Savart law as following:
\[
\mathbf{v}=\omega\,*\frac{1}{2\pi}\frac{-x^\bot}{|x|^2},
\]
where $x^\bot=(x_2,\,-x_1)$ if $x=(x_1,\,x_2)$. One special singular
solutions of Euler equations is given by
$\omega=\sum^m_{i=1}\kappa_i\delta_{x_i(t)}$, which is related
\[
\mathbf{v}=-\sum^m_{i=1}\frac{\kappa_i}{2\pi}\frac{(x-x_i(t))^\bot}{|x-x_i(t)|^2}.
\]
and the positions of the vortices $x_i: \mathbb{R}\rightarrow
\mathbb{R}^2$ satisfy the following Kirchhoff law:
\[
\kappa_i\,\frac{dx_i}{dt}=(\nabla_{x_i}\mathcal{W})^\bot
\]
where $\mathcal{W}$ is the so called Kirchhoff-Routh function
defined by
\[
\mathcal{W}(x_1,\cdots, x_m)=\frac{1}{2}\sum_{i\neq
j}^m\frac{\kappa_i\kappa_j}{2\pi}\log\frac{1}{|x_i-x_j|}.
\]

In simply-connected bounded domain $\Omega\subset \mathbb{R}^2$,
similar singular solutions also exist. Suppose that the normal
component of $\mathbf{v}$ vanishes on $\partial\Omega$, then the
Kirchhoff-Routh function is

\begin{equation}\label{W1}
\mathcal{W}(x_1,\cdots, x_m)=\frac{1}{2}\sum_{i\neq
j}^m{\kappa_i\kappa_j}G(x_i,\,x_j)+
\frac{1}{2}\sum_{i=1}^m{\kappa_i^2}H(x_i,\,x_i),
\end{equation}
where $G$ is the Green function of $-\Delta$ on $\Omega$ with 0
Dirichlet boundary condition and $H$ is its regular part (the Robin
function).

 Let $v_n$ be the outward component of the velocity $\mathbf{v}$ on
 the boundary $\partial\Omega$, then we see that $\int_{\partial\Omega}v_n=0$ due to
 the fact that $\nabla\cdot\mathbf{v}=0$. Suppose that
 $\mathbf{v}_0$ is the unique harmonic field whose normal component
 on the boundary $\partial\Omega$ is $v_n$, then $\mathbf{v}_0$
 satisfies

\begin{equation}\label{v0}
\begin{cases}
 \nabla\cdot\mathbf{v}_0=0,\,\,\text{in}\,\Omega,\\
\nabla\times\mathbf{v}_0=0,\,\,\text{in}\,\Omega,\\
n\cdot\mathbf{v}_0=v_n,\,\, \text{on}\,\partial\Omega.
\end{cases}
\end{equation}

If $\Omega$ is simply-connected, then $\mathbf{v}_0$ can be written
$\mathbf{v}_0=(\nabla\psi_0)^\bot$, where the stream function
$\psi_0$ is determined up to a constant by
\begin{equation}\label{psi}
\begin{cases}
 -\Delta \psi_0=0,\,\,\text{in}\,\Omega,\\
-\displaystyle\frac{\partial\psi_0}{\partial\tau}=v_n,\,\,
\text{on}\,\partial\Omega,
\end{cases}
\end{equation}
where $\frac{\partial\psi_0}{\partial\tau}$ denotes the tangential
derivative on $\partial\Omega$. The Kirchhoff-Routh function
associated to the vortex dynamics becomes(see Lin \cite{Lin})

\begin{equation}\label{KR0}
\mathcal{W}(x_1,\cdots,x_m)=\frac{1}{2}\sum_{i\neq
j}^m\kappa_i\kappa _jG(x_i,x_j)+\frac{1}{2}\sum^{m}_{i=1}\kappa^2
_iH(x_i,x_i)+\sum^{m}_{i=1}\kappa_i\psi_0(x_i).
\end{equation}

 For $m$ clockwise vortices motion (corresponding to $\kappa^+_i>0$) and $n$ anti-clockwise vortices motion (corresponding to
 $-\kappa^-_j<0$), the
Kirchhoff-Routh function associated to the vortex dynamics becomes
\begin{equation}\label{KR}
\begin{split}
\mathcal{W}(x^+_1,\cdots,x^+_m,x^-_1,\cdots,x^-_n)
 =&\frac{1}{2}\sum_{i,k=1,i\neq
k}^m\kappa_i^+\kappa
_k^+G(x_i^+,x_k^+)+\frac{1}{2}\sum^{n}_{j,l=1,j\neq
l}\kappa^-_j\kappa^-_l
G(x_j^-,x_l^-)\\
 &+\frac 12\sum_{i=1}^m(\kappa_i^+)^2H(x_i^+,x_i^+)+\frac
12\sum_{j=1}^n(\kappa_j^-)^2H(x_j^-,x_j^-)\\
&-\sum_{i=1}^m\sum_{j=1}^n\kappa_i^+\kappa_j^-G(x_i^+,x_j^-)
+\sum^{m}_{i=1}\kappa_i^+\psi_0(x_i^+)-\sum^{n}_{j=1}\kappa_j^-\psi_0(x_j^-).
\end{split}
\end{equation}

It is known that critical points of the Kirchhoff-Routh function
$\mathcal{W}$
 give rise to stationary vortex points solutions of the Euler
 equations. As for the existence of critical points of $\mathcal{W}$ given by \eqref{W1},
 we refer to \cite{BT}.

 Roughly speaking, there are two methods to construct stationary
solutions of the Euler equation, which are the vorticity method and
the stream-function method. The vorticity method was first
established by Arnold and Khesin \cite{AK} and further developed by
Burton \cite{B} and Turkington \cite{T}.

The stream-function method consists in observing that if $\psi$
satisfies  $-\Delta \psi=f(\psi)$ for some function $f\in
C^1(\mathbb{R})$, then $\mathbf{v}=(\nabla\psi)^\bot$ and
$P=F(\psi)-\frac{1 }{2}|\nabla\psi|^2$ is a stationary solution to
the Euler equations, where
$(\nabla\psi)^\bot:=(\frac{\partial\psi}{\partial
x_2},-\frac{\partial\psi}{\partial x_1}), F(t)=\int_0^tf(s)ds$.
Moreover, the velocity $\mathbf{v}$ is irrotational on the set where
$f(\psi)=0$.

Set $q=-\psi_0$ and $u=\psi-\psi_0$, then $u$ satisfies the
following boundary value problem
\begin{equation}\label{1.3}
\begin{cases}
-\Delta u=f(u-q),\quad & x\in\Omega,\\
u=0,\quad &x\in\partial\Omega.
\end{cases}
\end{equation}
 In addition, if we suppose that $\inf_\Omega q>0$ and
$f(t)=0,~t\leq0$, the vorticity set $\{x: f(\psi)>0\}$ is bounded
away from the boundary.

The motivation to study \eqref{1.3} is to justify the weak
formulation for point vortex solutions of the incompressible Euler
equations by approximating these solutions with classical solutions.

 Marchioro and Pulvirenti \cite{MP} have approximated these solutions
 on finite time intervals by considering regularized initial data
 for the vorticity. On the other hand, the stationary point vortex
 solutions can also be approximated by stationary classical
 solutions. See e.g. \cite{AS,AY,Ba,BF1,FB,Ni,N,SV,T,Y1,Y2} and the references therein.

In \cite{EM3} Elcrat and  Miller, by a rearrangements of functions,
have studied steady, inviscid flows in two dimensions which have
concentrated regions of vorticity. In particular, they studied such
flows which "desingularize" a configuration of point vortices in
stable equilibrium with an irrotational flow, which generalized
their earlier work for one vortex \cite{EM1}\cite{EM2} which in turn
were based on results of Turkington \cite{T}. As pointed by Elcrat
and  Miller, an essential hypothesis in their existence proof was
that the vorticity was in a neighborhood of a stable point vortex
configuration.  Saffman and Sheffield \cite{SS} have found an
example of a steady flow in aerodynamics with a single point vortex
which is stable for a certain range of the parameters. This has been
generalized in \cite{EM1}, where some examples computationally of
stable configurations of two point vortices were briefly discussed.
Further examples of multiple point vortex configurations are given
in \cite{M1}, where a theorem on the existence of such
configurations is also given.

  It is worth
 pointing out that except \cite{EM3} the above approximations can just give
 explanation for the formulation to single point vortex solutions.
D. Smets and J. Van Schaftingen \cite{SV} investigated the following
problem
\begin{equation}\label{00}
\begin{cases}
-\ep^2 \Delta u=\left(u-q-\frac{\kappa}{2\pi}\ln\frac{1}{\ep}\right)_+^{p},  & \text{in}\;\Om,\\
u=0, &\text{on}\; \partial\Om,
\end{cases}
\end{equation}
and gave the exact asymptotic behavior and expansion  of the least
energy solution by estimating the upper bounds on the energy. The
solutions for \eqref{00} in \cite{SV} were obtained by finding a
minimizer of the corresponding functional in a suitable function
space, which can only give approximation to a single point
non-vanishing vortex.  In \cite{CLW}, we have shown that multi-point
vortex solutions can
 be approximated by stationary classical solutions.

Concerning regularization of pairs of vortices, D. Smets and J. Van
Schaftingen \cite{SV} also studied the following problem
 \begin{equation}\label{11}
\begin{cases}
-\ep^2 \Delta u=\left(u-q-\frac{\kappa}{2\pi}\ln\frac{1}{\ep}\right)_+^{p}-(q-\frac{\kappa}{2\pi}\ln\frac{1}{\ep}-u)_+^p,  & \text{in}\;\Om,\\
u=0, &\text{on}\; \partial\Om,
\end{cases}
\end{equation}
and obtained the exact asymptotic behavior and expansion  of the
least energy solution by similar methods for \eqref{00}. This method
is hard  to obtain multiple non-vanishing pairs of vortices
solutions.

 In this paper, we approximate stationary vortex solutions of Euler equations \eqref{1.2} with multiple non-vanishing pairs of vortices
solutions
 by stationary classical solutions. Our main result concerning \eqref{1.2} is the
 following:

\begin{theorem}\label{nth1}

Suppose that $\Omega\subset \mathbb{R}^2$ is a bounded
simply-connected smooth domain. Let $v_n:
\partial\Omega\rightarrow \mathbb{R}$ be such that $v_n\in
L^s(\partial\Omega)$ for some $s>1$ satisfying
$\int_{\partial\Omega}v_n=0$. Let
$\kappa^+_i>0,\kappa^-_j>0,~i=1,\cdots,m,~j=1,\dots,n$. Then, for
any given $C^1$-stable critical point
$(x_{1,\ast}^+,\cdots,x_{m,\ast}^+,x_{1,\ast}^-,\cdots,x_{n,\ast}^-)$
of Kirchhoff-Routh function $\mathcal{W}$ defined by \eqref{KR},
there exists $\ep_0>0$, such that for each $\ep\in (0,\ep_0)$,
problem \eqref{1.2} has a stationary solution $\mathbf{v}_\ep$ with
outward boundary flux given by $v_n$, such that its vorticities
$\omega^\pm_\ep$ satisfying
\[
supp (\omega^+_\ep)\subset\cup_{i=1}^m
B(x_{i,\,\ep},C\ep)~~\text{for}~~x^+_{i,\,\ep}\in\Omega,
~~i=1,\cdots,m,
\]
\[
supp (\omega^-_\ep)\subset\cup_{j=1}^n
B(x^-_{j,\,\ep},C\ep)~~\text{for}~~x^-_{j,\,\ep}\in\Omega,
~~j=1,\cdots,n
\]
and as $\ep\rightarrow0$
\[
\int_\Omega \omega_\ep\rightarrow
\sum_{i=1}^m\kappa_i^+-\sum_{j=1}^n\kappa^-_j,
\]

\[
(x^+_{1,\,\ep},\cdots,x^+_{m,\,\ep},x_{1,\,\ep}^-,\cdots,x^-_{n,\ep})\rightarrow
(x_{1,\ast}^+,\cdots,x_{m,\ast}^+,x_{1,\ast}^-,\cdots,x_{n,\ast}^-).
\]
\end{theorem}

\begin{remark} The simplest case, corresponding
 to pairs of vortices $(m=n=1)$  was studied by Smets and Van Schaftingen
 \cite{SV} by minimizing the corresponding energy functional. In their paper
as $\ep\to 0$, $\mathcal{W}(x^+_{1,\,\ep},x^-_{1,\ep})\rightarrow
 \sup\limits_{x_1^+,x^-_1\in\Omega,\, x_1^+\neq x^-_1}\mathcal{W}(x_1^+,x_1^-)$. Even in the case $m=n=1$, our result
 extends theirs to general critical points (with additional assumption that
 the critical point is
 non-degenerate or stable in the sense of $C^1$). The method used in \cite{SV} can not be applied to deal with
 general critical point cases. The method used here is constructive and  is completely
 different from theirs.
\end{remark}

\begin{remark}
In this case that $m=n=1$ suppose that
$(x^+_{1,\ast},\,x^-_{1,\ast})$ is a strict local maximum(or
minimum) point of Kirchhoff-Routh function $\mathcal{W}(x^+,\,x^-)$
defined by \eqref{KR}, statement of Theorem \ref{nth1} still holds
which can be proved similarly (see Remark \ref{re1.4}). Thus we can
obtain corresponding existence result in \cite{SV}.
\end{remark}

Theorem \ref{nth1} is proved via considering the following  problem
\begin{equation}\label{0}
 \begin{cases}
-\ep^2 \Delta
u=\sum_{i=1}^m\chi_{\Om_i^+}\big(u-q-\frac{\kappa_i^+}{2\pi}\ln\frac{1}{\ep}\big)_+^p
-\sum_{j=1}^n\chi_{\Om_j^-}\big(q-\frac{\kappa_j^-}{2\pi}\ln\frac{1}{\ep}-u\big)_+^p, \quad & x\in\Omega, \\
u=0, \quad & x\in\partial\Omega,
\end{cases}
\end{equation}
where $p>1$, $q\in C^2(\Omega)$,  $\Omega\subset\mathbb{R}^2$ is a
bounded domain, $\Om_i^+(i=1\cdots,m)$ and $\Om_j^-(\,j=1\cdots,n)$
are mutually disjoint subdomains of $\Omega$ such that
$x^+_{i,*}\in\Om^+_i,$ and $x^-_{j,*}\in\Om^-_j$.

\begin{theorem}\label{th3}

Suppose $q\in C^2(\Omega)$. Then for any given
$\kappa^+_i>0,\kappa^-_j>0,~i=1,\cdots,m,~j=1,\dots,n$ and for any
given $C^1$-stable critical point
$(x_{1,\ast}^+,\cdots,x_{m,\ast}^+,x_{1,\ast}^-,\cdots,x_{n,\ast}^-)$
of Kirchhoff-Routh function $\mathcal{W}$ defined by \eqref{KR},
there exists $\ep_0>0$, such that for each $\ep\in (0,\ep_0)$,
\eqref{0} has a solution $u_\ep$, such that the set
$\Omega_{\ep,i}^+ =\{x:
u_\ep(x)-\frac{\kappa_i^+}{2\pi}\,\ln\frac{1}{\ep}
-q(x)>0\}\subset\subset\Om_i^+, \,i=1,\cdots,m$, $\Omega_{\ep,j}^-
=\{x: u_\ep(x)-\frac{\kappa_j^-}{2\pi}\,\ln\frac{1}{\ep}
-q(x)>0\}\subset\subset\Om_j^-, \,j=1,\cdots,n$ and as $\ep\to 0$,
each $\Omega^\pm_{\ep,\,i}$ shrinks to $x_{i,\ast}^\pm\in \Omega$.

\end{theorem}

\begin{remark}\label{re1.4}
For the case $m=n=1$, suppose that $(x^+_{1,\ast},\,x^-_{1,\ast})$
is a strict local maximum(or minimum) point of Kirchhoff-Routh
function $\mathcal{W}(x)$ defined by \eqref{KR}, then statement of
Theorem \ref{th3} still holds. This conclusion can be proved by
making corresponding modification of the proof of Theorem \ref{th3}
in obtaining critical point of $K(z)$ defined by \eqref{K}(see
Propositions 2.3,2.5 and 2.6 in \cite{Cao} for detailed arguments).
\end{remark}

As in \cite{CLW},
 we prove  Theorem \ref{th3}  by considering an equivalent problem of \eqref{0}
 instead. Let $ w=\frac{2\pi}{|\ln\ep|}u$ and
$\delta=\ep(\frac{2\pi}{ |\ln\ep |})^{\frac{p-1}{2}}$, then
\eqref{0} becomes
\begin{equation}\label{1}
\begin{cases}
-\de^2 \Delta w=\sum_{i=1}^m\chi_{\Om_i^+}\left(w-\kappa_i^+-\frac{2\pi}{|\ln\ep|}q(x)\right)_+^{p}
-\sum_{j=1}^m\chi_{\Om_j^-}\left(\frac{2\pi}{ |\ln\ep|}q(x)-\kappa_j^--w\right)_+^{p},  & \text{in}\;\Om,\\
w=0, &\text{on}\; \partial\Om.
\end{cases}
\end{equation}

We will use a reduction argument to prove Theorem~\ref{th3}. To
 this end, we need to construct an approximate solution for
\eqref{1}.  For the problem studied in this paper, the corresponding
``limit" problem in $\mathbb{R}^2$ has no bounded nontrivial
solution. So, we will follow the method in \cite{CPY,DY} to
construct an approximate solution. Since there are two parameters
$\de,~\ep$ in problem \eqref{1} and two terms in nonlinearity, which
causes some difficulty, we must take this influence into careful
consideration and give delicate estimates in order to perform the
reduction argument. For example we need to consider
$(s^+_{1,\de},\cdots,s^+_{m,\de},s^-_{1,\de},\cdots,s^-_{n,\de})$
and
$(a^+_{1,\de},\cdots,a^+_{m,\de},a_{1,\de}^-,\cdots,a_{n,\de}^-)$
together in Lemma~\ref{l2.1}.

As a final remark, we point out that problem \eqref{1} can  be
considered as a free boundary problem. Similar problems have been
studied extensively. The reader can refer to
\cite{CF,CLW,CPY,DY,FW,LP} for more results on this kind of
problems.

This paper is organized as follows.  In section~2, we construct the
approximate solution for \eqref{1}. We will carry out a reduction
argument in section~3 and the main results will be proved in
section~4.  We put some basic estimates used in sections 3 and 4 in
the appendix.

\section{Approximate solutions }

In the section,  we will construct approximate solutions for
\eqref{1}.

Let $R>0$ be a large constant, such that for any $x\in \Om$,
$\Om\subset\subset B_R(x)$. Consider the following problem:

\begin{equation}\label{2.1}
\begin{cases}
-\delta^2\Delta w=( w-a)_+^{p},& \text{in}\; B_R(0),\\
w=0, &\text{on}\;\partial B_R(0),
\end{cases}
\end{equation}
where $a>0$ is a constant. Then, \eqref{2.1} has a unique solution
$W_{\de,a}$, which can be written as

\begin{equation}\label{2.2}
W_{\de,a}(x)=
\begin{cases}
a+\de^{2/(p-1)}s_\de^{-2/(p-1)}\phi\bigl(\frac{|x|}{s_\de}\bigr), &  |x|\le s_\de,\\
a\ln\frac {|x|} R/\ln \frac {s_\de}R, & s_\de\le |x|\le R,
\end{cases}
\end{equation}
where $\phi(x)=\phi(|x|)$ is the unique solution of
\begin{equation*}
-\Delta \phi=\phi^{p},\quad\phi>0,~~\phi\in H_0^1\bigl(B_1(0)\bigr)
\end{equation*}
and $s_\de\in (0,R)$ satisfies
\[
\de^{2/(p-1)}s_\de^{-2/(p-1)}\phi^\prime(1)=\frac{a}{\ln(s_\de/R)},
\]
which implies
\[
\frac{s_\de}{\de|\ln\de|^{(p-1)/2}}\rightarrow\left(\frac{|\phi^\prime(1)|}{a}\right)^{(p-1)/2}>0,\quad\text{as}~~\de\rightarrow0.
\]
Moreover, by Pohozaev identity, we can get that
\[
\int_{B_1(0)}\phi^{p+1}=\frac{\pi(p+1)
}{2}|\phi^\prime(1)|^2~~\text{and}~~\int_{B_1(0)}\phi^{p}=2\pi|\phi^\prime(1)|.
\]

 For any $z \in \Om$,  define
$W_{\de,z,a}(x)=W_{\de,a}(x-z)$. Because $W_{\de,z,a}$
 does not vanish on $\partial\Omega$, we need to  make a  projection. Let
$PW_{\de,z,a}$ be the solution of

\[
\begin{cases}
-\de^2 \Delta w=( W_{\de,z,a}-a)_+^{p},& \text{in } \; \Om,\\
w=0, &\text{on}\; \partial \Om.
\end{cases}
\]
Then

\begin{equation}\label{2.3}
PW_{\de,z,a}= W_{\de,z,a}-\frac a{\ln \frac{R}{s_\de}} g(x,z),
\end{equation}
where $g(x,z)$ satisfies

\[
\begin{cases}
- \Delta g=0,& \text{in } \; \Om,\\
g=\ln\frac{R}{|x-z|}, &\text{on}\; \partial \Om.
\end{cases}
\]
It is easy to see that

\[
g(x,z)=\ln R +2\pi h(x,z),
\]
where $h(x,z)=-H(x,z)$.

Let $Z=(Z_m^+,Z_n^-)$, where $Z_m^+=(z_1^+,\cdots,z_m^+)$,
$Z_n^+=(z_1^-,\cdots,z_n^-)$. We will construct solutions for
\eqref{1} of  the form

\[
\sum_{i=1}^m PW_{\de,z_{i}^+,a_{\de,i}^+}-\sum_{j=1}^n
PW_{\de,z_{j}^-,a_{\de,j}^-} +\omega_\de,
\]
where $z_i^+,z_j^-\in\Omega$, $a^+_{\de,i}>0,a^-_{\de,j}>0$ for
$i=1,\cdots,m$, $j=1,\cdots,n$, $\omega_\de$ is a perturbation term.
To make $\omega_\de$ as small as possible, we need to choose
$a^+_{\de,i},\,a^-_{\de,j}$ properly.

In this paper, we always assume that $z_i^+, z_j^-\in\Om$ satisfies

\begin{equation}\label{2.5}
\begin{split}
&d(z_i^+,\partial\Om)\ge \varrho,~d(z_j^-,\partial\Om)\ge
\varrho,\quad |z_i^+-z_k^+|\ge \varrho^{\bar L},\quad i, k=1,\cdots,m,\; i\ne k\\
&~~|z_j^--z_l^-|\ge \varrho^{\bar L},\quad|z_i^+-z_j^-|\ge
\varrho^{\bar L},\quad j, l=1,\cdots,n,\; j\ne l,
\end{split}
\end{equation}
where $\varrho>0$ is a fixed small constant and $\bar L>0$ is a
fixed large constant.

\begin{lemma}\label{l2.1}
For $\de>0$ small,  there exist
$(s_{\de,1}^+(Z),\cdots,s_{\de,m}^+(Z),s_{\de,1}^-(Z),\cdots,s^-_{\de,n}(Z))$
and
$(a_{\de,1}^+(Z),\cdots,a_{\de,m}^+(Z),a_{\de,1}^-(Z),\cdots,a^-_{\de,n}(Z))$
satisfying the following system
\begin{equation}\label{2.6}
\de^{2/(p-1)}(s_i^+)^{-2/(p-1)}\phi^\prime(1)=\frac{a_{i}^+}{\ln(s_i^+/R)},\qquad
i=1,\cdots,m
\end{equation}
\begin{equation}\label{2.61}
\de^{2/(p-1)}(s_j^-)^{-2/(p-1)}\phi^\prime(1)=\frac{a_{j}^-}{\ln(s_j^-/R)},\qquad
j=1,\cdots,n
\end{equation}
and
\begin{equation}\label{2.7}
a_{i}^+= \kappa_i^+ + \frac{2\pi
q(z_i^+)}{|\ln\ep|}+\frac{g(z_i^+,z_i^+)}{\ln\frac R {s_{i}^+} }
a_{i}^+ -\sum_{\alpha\ne i}^m \frac{\bar
G(z_i^+,z_\alpha^+)}{\ln\frac R{s_{\alpha}^+}
}a_{\alpha}^++\sum_{l=1}^n\frac{\bar
G(z_i^+,z_l^-)}{\ln\frac{R}{s_l^-}}a_l^- ,\quad i=1,\cdots, m,
\end{equation}
\begin{equation}\label{2.71}
a_{j}^-= \kappa_j^- - \frac{2\pi
q(z_j^-)}{|\ln\ep|}+\frac{g(z_i^-,z_i^-)}{\ln\frac R {s_{j}^-} }
a_{j}^- -\sum_{\beta\ne j}^n \frac{\bar G(z_\beta^-,z_j^-)}{\ln\frac
R{s_{\beta}^+} }a_{\beta}^++\sum_{k=1}^m\frac{\bar
G(z_j^-,z_k^+)}{\ln\frac{R}{s_k^-}}a_k^- ,\quad j=1,\cdots, n,
\end{equation}
where $ \bar G(x,y)=\ln \frac R{|x-y|} -g(x,y)$ for $x\neq y$.
\end{lemma}

Since the proof is exactly the same as in Lemma 2.1 in \cite{CLW},
we omit it here therefore.

To simplify our notations, for given $Z=(Z_m^+,Z_n^-)$, in this
paper, we will use $a_{\de,i}^\pm$, $s_{\de,i}^\pm$ instead of
$a_{\de,i}^\pm(Z)$, $s_{\de,i}^\pm(Z)$.  From now on we will always
choose
$(a_{\de,1}^+,\cdots,a_{\de,m}^+,a_{\de,1}^-,\cdots,a_{\de,n}^-)$
and
$(s_{\de,1}^+,\cdots,s_{\de,m}^+,s_{\de,1}^-,\cdots,s_{\de,n}^-)$
such that \eqref{2.6}--\eqref{2.71} hold. For
$(a_{\de,1}^+,\cdots,a_{\de,m}^+,a_{\de,1}^-,\cdots,a_{\de,n}^-)$
and
$(s_{\de,1}^+,\cdots,s_{\de,m}^+,s_{\de,1}^-,\cdots,s_{\de,n}^-)$
chosen in such a way let us define
\begin{equation}\label{2.8}
P_{\de, Z,i}^+= PW_{\de,z^+_i,\,a_{\de,i}^+},~~P_{\de, Z,j}^-=
PW_{\de,z^-_j,\,a_{\de,j}^-}.
\end{equation}

\begin{remark}\label{remark2.2}
As in \cite{CLW}, we have the following asymptotic expansions:
\begin{equation}\label{r2.2.1}
\frac{1}{\ln\frac{R}{s^{+}_{\de,i}}}=\frac{1}{\ln\frac{R}{\ep}}
+0\left(\frac{\ln|\ln\ep|}{|\ln\ep|^2}\right),\,i=1,\cdots,m,
\end{equation}

\begin{equation}\label{r2.2.2}
a_{\de,i}^+=1+\frac{2\pi
q(z_i^+)}{\kappa|\ln\ep|}+\frac{g(z_i^+,z_i^+)}{\ln\frac{R}{\ep}}-\sum_{\alpha\neq
i}^m\frac {\bar
G(z_i^+,z_\alpha^+)}{\ln\frac{R}{\ep}}+\sum_{l=1}^n\frac{\bar
G(z_i^+,z_l^-)}{\ln\frac{R}{\ep}}+
0\Bigl(\frac{\ln|\ln\ep|}{|\ln\ep|^2}\Bigr),\, i=1,\cdots,m,
\end{equation}

\begin{equation}\label{2.10}
\frac{\partial a_{\de,i}^+}{\partial
z_{k,h}^\pm}=0\left(\frac1{|\ln\ep|}\right),\quad~~\frac{\partial
s_{\de,i}^+}{\partial
z_{k,h}^\pm}=0\left(\frac{\ep}{|\ln\ep|}\right),\,i=1,\cdots,m,\,h=1,2.
\end{equation}
Moreover, $a_{\de,j}^-$ and $s_{\de,j}^-$ have similar expansions.

\end{remark}

To simplify notations, set
\[
P^+_{\de,Z}=\sum_{\alpha=1}^m P^+_{\de,
Z,\alpha},~~P^-_{\de,Z}=\sum_{\beta=1}^n P^-_{\de,Z,\beta}.
\]

 Then,  we find that for $x\in
B_{L s_{\de,i}^+}(z_i^+)$, where $L>0$ is any fixed constant,

\[
\begin{split}
& P_{\de, Z,i}^+(x)-\kappa_i^+-\frac{2\pi q(x)}{|\ln\ep|}=
W_{\de,z^+_i,\, a_{\de,i}^+}(x) -\frac {a^+_{\de,i}}{\ln
\frac{R}{s^+_{\de,i} }} g(x,z_i^+)-\kappa_i^+-\frac{2\pi
q(x)}{|\ln\ep|}
\\
=&  W_{\de,z^+_i,\, a_{\de,i}^+}(x)-\kappa_i^+-\frac
{a^+_{\de,i}}{\ln \frac{R}{s^+_{ \de,i}}}
 g(z_i^+,z_i^+)-\frac {a^+_{\de,i}}{\ln \frac{R}{s^+_{\de,i}}}\Bigl(
 \left\langle D g(z_i^+,z_i^+), x-z_i^+\right\rangle +O( |x-z_i^+|^2)\Bigr)\\
\quad&-\frac{2\pi q(z_i^+)}{|\ln\ep|}-\frac{2\pi}{|\ln\ep|}\left(\left\langle Dq(z_i^+),x-z_i^+\right\rangle+O(|x-z_i^+|^2)\right)\\
 =&  W_{\de,z^+_i,\, a^+_{\de,i}}(x)- \kappa_i^+-\frac{2\pi q(z^+_i)}{|\ln\ep|}-\frac{2\pi}{|\ln\ep|}\left\langle
Dq(z_i^+),x-z^+_i\right\rangle\\
\quad&-\frac {a^+_{\de,i}}{\ln \frac{R}{s^+_{ \de,i}}}
 g(z_i^+,z_i^+)-\frac {a^+_{\de,i}}{\ln \frac{R}{s^+_{\de,i}}}
 \left\langle D g(z_i^+,z_i^+), x-z_i^+\right\rangle +O\left(\frac{(s^+_{\de,i})^2}{|\ln\ep|}\right),
\end{split}
\]
and for $k\ne i$ and  $x\in B_{Ls^+_{\de,i}}(z_i^+)$, by \eqref{2.2}

\[
\begin{split}
&  P^+_{\de, Z,k}(x)=W_{\de,
z_k^+,a^+_{\de,k}}(x)-\frac{a^+_{\de,k}}{\ln\frac R{s^+_{\de,k}} }
  g(x,z_k^+)=
 \frac{a^+_{\de,k}}{\ln\frac R{s^+_{\de,k}}
}   \bar G(x,z_k^+)\\
=&
  \frac{a^+_{\de,k}}{\ln\frac R{s^+_{\de,k}}
}     \bar  G(z_i^+,z_k^+)+\frac{a^+_{\de,k}}{\ln\frac
R{s^+_{\de,k}} } \left\langle D  \bar
G(z_i^+,z_k^+),x-z_i^+\right\rangle + O\Bigl( \frac{(s^+_{\de,i})^2
}{|\ln\ep| } \Bigr)
\end{split}
\]
and

\[
 P^-_{\de, Z,j}(x)=\frac{a^-_{\de,j}}{\ln\frac R{s^-_{\de,j}}
}     \bar  G(z_i^+,z_j^-)+\frac{a^-_{\de,j}}{\ln\frac
R{s^-_{\de,j}} } \left\langle D  \bar
G(z_i^+,z_j^-),x-z_i^+\right\rangle +
 O\Bigl( \frac{(s^+_{\de,i})^2 }{|\ln\ep|
} \Bigr).
\]
 So,  by using \eqref{2.7}, we obtain

\begin{equation}\label{2.9}
\begin{split}
&  P^+_{\de, Z}(x)-P^-_{\de,Z}(x)-\kappa_i^+-\frac{2\pi q(x)}{|\ln\ep|}\\
=&
 W_{\de,z_i^+, a^+_{\de,i}}(x)- a^+_{\de,i} -\frac{2\pi}{|\ln\ep|}\left\langle Dq(z_i^+),x-z^+_i\right\rangle-\frac {a^+_{\de,i}}{\ln \frac{R}{s^+_{\de,i}}}
 \left\langle D g(z_i^+,z_i^+), x-z^+_i\right\rangle \\
& +\sum_{k\ne i}^m \frac{a^+_{\de,k} }{\ln\frac R{s^+_{\de,k}} }
\left\langle D \bar G(z_i^+,z_k^+),x-z_i^+\right\rangle-\sum_{l=1}^n
\frac{a^-_{\de,l} }{\ln\frac R{s^-_{\de,l}} } \left\langle D \bar
G(z_i^+,z_l^-),x-z_i^+\right\rangle\\
&  +O\left( \frac{(s^+_{\de,i})^2 }{|\ln\ep| } \right),\quad x\in
B_{L s^+_{\de,i}}(z^+_i).
\end{split}
\end{equation}
Similarly, we have
\begin{equation}\label{2.91}
\begin{split}
& P^-_{\de,Z}(x)- P^+_{\de, Z}(x)-\kappa_j^-+\frac{2\pi q(x)}{|\ln\ep|}\\
=&
 W_{\de,z_j^-, a^-_{\de,j}}(x)- a^-_{\de,j} +\frac{2\pi}{|\ln\ep|}\left\langle Dq(z_j^-),x-z^-_j\right\rangle-\frac {a^-_{\de,j}}{\ln \frac{R}{s^-_{\de,j}}}
 \left\langle D g(z_j^-,z_j^-), x-z^-_j\right\rangle \\
& +\sum_{l\ne j}^n \frac{a^-_{\de,l} }{\ln\frac R{s^-_{\de,l}} }
\left\langle D \bar G(z_j^-,z_l^-),x-z_j^-\right\rangle-\sum_{k=1}^m
\frac{a^+_{\de,k} }{\ln\frac R{s^+_{\de,k}} } \left\langle D \bar
G(z_j^-,z_k^+),x-z_j^-\right\rangle\\
&  +O\left( \frac{(s^-_{\de,j})^2 }{|\ln\ep| } \right),\quad x\in
B_{L s^-_{\de,j}}(z_j^-).
\end{split}
\end{equation}

We end this section by giving the following formula which can be
obtained by direct computation and will be used in the next two
sections.

\begin{equation}\label{2.11}
\begin{array}{ll}
 \displaystyle\frac{\partial W_{\de,z_i^\pm,a^\pm_{\de,i}}(x)}{\partial z^\pm_{i,h}}&\\
 =\left\{
 \begin{array}{ll}
\displaystyle\frac{1}{\de}\Bigl(\frac{a^\pm_{\de,i}}{|\phi^\prime(1)||\ln\frac{R}{s^\pm_{\de,i}}|}\Bigr)^{(p+1)/2}
\phi^\prime\bigl(\frac{|x-z^\pm_i|}{s^\pm_{\de,i}}\bigr)
\frac{z^\pm_{i,h}-x_h}{|x-z_i^\pm|}+
O\left(\frac{1}{|\ln\ep|}\right),
 ~~x\in B_{s^\pm_{\de,i}}(z^\pm_i),\\
 \,\\
\displaystyle-\frac{a^\pm_{\de,i}}{\ln\frac{R}{s^\pm_{\de,i}}}\frac{z^\pm_{i,h}-x_h}{|x-z^\pm_i|^2}+O\left(\frac{1}{|\ln\ep|}\right),
 \qquad\qquad\qquad\qquad\qquad\quad x\in \Omega\setminus B_{s^\pm_{\de,i}}(z^\pm_i).
\end{array}
\right.\\
\end{array}
\end{equation}

\section{the reduction}

Let

\[
w(x)=
\begin{cases}
\phi(|x|), &|x|\le 1,\\
\phi^\prime(1)\ln |x|,  & |x|>1.
\end{cases}
\]
Then $w\in C^1(\mathbb{R}^2)$. Since $\phi^\prime(1)<0$ and $\ln
|x|$ is harmonic for $|x|>1$, we see that $w$ satisfies

\begin{equation}\label{3.1}
-\Delta w= w_+^{p}, \quad \text{in}\; \mathbb{R}^2.
\end{equation}
Moreover, since $w_+$ is Lip-continuous,  by the Schauder estimate,
$w\in C^{2,\alpha}$ for any $\alpha\in (0,1)$.

Consider the following problem:

\begin{equation}\label{3.2}
-\Delta v-  pw_+^{p-1} v=0,\quad v\in L^\infty(\mathbb{R}^2),
\end{equation}
It is easy to see  that $\frac{\partial w}{\partial x_i}$, $i=1,2,$
is a solution of \eqref{3.2}. Moreover, from Dancer and Yan
\cite{DY}, we know that $w$ is also non-degenerate, in the sense
that the kernel of the operator $Lv:=-\Delta v- pw_+^{p-1} v,~~v\in
D^{1,2}(\mathbb{R}^2)$ is spanned by $\bigl\{\frac{\partial
w}{\partial x_1}, \frac{\partial w}{\partial x_2}\bigr\}$.

Let $P_{\de,Z,i}^+,~P_{\de,Z,j}^-$ be the functions defined in
\eqref{2.8}. Set

\[
\begin{split}
F_{\de,Z}=\Bigg\{u: u\in L^p(\Om), & \int_\Om
 \frac{\partial
P^+_{\de,Z,i}}{\partial z^+_{i,h}}u =0, \int_\Om
 \frac{\partial
P^-_{\de,Z,j}}{\partial z^-_{j,h}}u =0,\\
& i=1,\cdots,m, \;\;j=1,\cdots,n,\;\;h=1,2 \Bigg\},
\end{split}
\]
and

\[
\begin{split}
E_{\de,Z}=\Bigg\{u:\;u\in W^{2,p}(\Om)\cap H_0^1(\Om), &\int_\Om
 \Delta\left( \frac{\partial
P^+_{\de,Z,i}}{\partial z^+_{i,h}}\right)u =0, \;\int_\Om
 \Delta\left( \frac{\partial
P^-_{\de,Z,j}}{\partial z^-_{j,h}}\right)u =0,\\
&i=1,\cdots,m,\;\;j=1,\cdots,n, \;\; h=1,2 \Bigg\}.
\end{split}
\]

For any $u\in L^p(\Om)$, define $Q_\de u$  as follows:
\[
Q_\de u= u-\sum_{i=1}^m \sum_{h=1}^2
b^+_{i,h}\left(-\de^2\Delta\Bigl(\frac{\partial
P^+_{\de,Z,i}}{\partial z^+_{i,h}}\Bigr)\right)-\sum_{j=1}^n
\sum_{\bar h=1}^2 b^-_{j,\bar
h}\left(-\de^2\Delta\Bigl(\frac{\partial P^-_{\de,Z,j}}{\partial
z^-_{j,\bar h}}\Bigr)
 \right),
\]
where  the constants $b^+_{i,h}$, $b^-_{j,\bar h}$  satisfy

\begin{equation}\label{3.4}
\begin{split}
&\sum_{i=1}^m \sum_{h=1}^2
b^+_{i,h}\left(-\de^2\int_\Om\Delta\Bigl(\frac{\partial
P^+_{\de,Z,i}}{\partial z^+_{i,h}}\Bigr)\frac{\partial
P^+_{\de,Z,k}}{\partial z^+_{k,\hat h}}\right)\\
&+\sum_{j=1}^n \sum_{\bar h=1}^2 b^-_{j,\bar
h}\left(-\de^2\int_\Om\Delta\Bigl(\frac{\partial
P^-_{\de,Z,j}}{\partial z^-_{j,\bar h}}\Bigr)\frac{\partial
P^+_{\de,Z,k}}{\partial z^+_{k,\hat h}}
 \right)=\int_\Om u\frac{\partial
P^+_{\de,Z,k}}{\partial z^+_{k,\hat h}},
 \end{split}
\end{equation}
and
\begin{equation}\label{3.41}
\begin{split}
&\sum_{i=1}^m \sum_{h=1}^2
b^+_{i,h}\left(-\de^2\int_\Om\Delta\Bigl(\frac{\partial
P^+_{\de,Z,i}}{\partial z^+_{i,h}}\Bigr)\frac{\partial
P^-_{\de,Z,l}}{\partial z^-_{l,\tilde h}}\right)\\
&+\sum_{j=1}^n \sum_{\bar h=1}^2 b^-_{j,\bar
h}\left(-\de^2\int_\Om\Delta\Bigl(\frac{\partial
P^-_{\de,Z,j}}{\partial z^-_{j,\bar h}}\Bigr)\frac{\partial
P^-_{\de,Z,l}}{\partial z^-_{l,\tilde h}}
 \right)=\int_\Om u\frac{\partial
P^-_{\de,Z,l}}{\partial z^-_{l,\tilde h}}.
 \end{split}
\end{equation}

Since  $\int_\Om
 \frac{\partial
P^+_{\de,Z,k}}{\partial z^+_{k,\hat h}} Q_\de u=0$, $\int_\Om
 \frac{\partial
P^-_{\de,Z,l}}{\partial z^-_{l,\tilde h}} Q_\de u=0$,  the operator
$Q_\de$ can be regarded as a projection from $L^p(\Omega)$ to
$F_{\de,Z}$. In order to show that we can solve \eqref{3.4} and
\eqref{3.41} to obtain $b^+_{i,h}$ and $b^-_{j,\bar h}$, we just
need the following estimate ( by \eqref{2.10} and \eqref{2.11}):

\begin{equation}\label{3.5}
\begin{split}
&-\de^2\int_{\Om}\Delta\Bigl( \frac{\partial P^+_{\de,Z,i}}{\partial
z^+_{i, h}}\Bigr)
\frac{\partial P^+_{\de,Z,k}}{\partial z^+_{k,\hat h}}\\
=&
p\int_{\Om}\bigl(W_{\de,z_i^+,a^+_{\de,i}}-a_{\de,i}^+\bigr)_+^{p-1}\left(\frac{\partial
W_{\de,z^+_{i},a^+_{\de,i}}}{\partial z^+_{i, h}}-\frac{\partial
a^+_{\de,i}}{\partial z^+_{i,h}}\right) \frac{\partial
P^+_{\de,Z,k}}{\partial z^+_{k,\hat h}}
\\
=& \delta_{ik h\hat h}\frac
c{|\ln\ep|^{p+1}}+0\left(\frac{\ep}{|\ln\ep|^{p+1}}\right),
\end{split}
\end{equation}
where $c>0$ is a constant, $\delta_{ikh \hat h}=1$,  if $i=k$ and
$h=\hat h$;  otherwise, $\delta_{ijh \hat h}=0$.

Similarly,
\begin{equation}\label{3.51}
\begin{split}
-\de^2\int_{\Om}\Delta\Bigl( \frac{\partial P^-_{\de,Z,j}}{\partial
z^-_{j, \bar h}}\Bigr) \frac{\partial P^-_{\de,Z,l}}{\partial
z^-_{l,\tilde h}} = \delta_{jl \bar h\tilde h}\frac
c{|\ln\ep|^{p+1}}+0\left(\frac{\ep}{|\ln\ep|^{p+1}}\right),
\end{split}
\end{equation}
where $c>0$ is a constant, $\delta_{jl\bar h\tilde h}=1$,  if $j=l$
and $\bar h=\tilde h$;  otherwise, $\delta_{jl\bar h \tilde h}=0$.

 Set

\[
\begin{split}
L_\de u= -\de^2 \Delta u&-\sum_{i=1}^m p\chi_{\Om_i^+}\left(
P^+_{\de,Z}-P^-_{\de,Z}-\kappa_i^+-\frac{2\pi
q(x)}{|\ln\ep|}\right)_+^{p-1}u\\
&-\sum_{j=1}^m p\chi_{\Om_j^-}\left(P_{\de,Z}^-
-P_{\de,Z}^+-\kappa_j^-+\frac{2\pi q(x)}{|\ln\ep|}\right)_+^{p-1}u,
\end{split}
\]
and
\[B_{\de,Z}=\Bigl(\cup_{i=1}^mB_{Ls_{\de,i}^+}(z_i^+)\Bigr)\bigcup\Bigl(\cup_{j=1}^nB_{Ls_{\de,j}^-}(z_j^-)\Bigr).\]

  We have the following lemma.

\begin{lemma}\label{l31}
 There are constants $\rho_0>0$ and $\de_0>0$, such that for any
$\de\in (0,\de_0]$,
 $Z$ satisfying \eqref{2.5},  $ u\in E_{\de,Z}$ with
$Q_\de L_\de u =0$ in $\Omega\setminus B_{\de,Z}$ for some $L>0$
large, then

\[
\|Q_\de  L_\de u\|_{L^p(\Om)}  \ge
\frac{\rho_0\de^{\frac{2}{p}}}{|\ln\de|^{\frac{(p-1)^2}{p}}}
\|u\|_{L^\infty(\Om)}.
\]

\end{lemma}

\begin{proof}

 Set $s_{N,j}^{\pm}=s_{\de_N,j}^{\pm}$.
In the sequel, we will use $\|\cdot\|_p, \|\cdot\|_\infty$ to denote
$\|\cdot\|_{L^p(\Om)}$ and  $ \|\cdot\|_{L^\infty(\Om)}$
respectively.

We argue by contradiction. Suppose that there are $\de_N\to 0$,
$Z_N$ satisfying \eqref{2.5} and $u_N\in E_{\de_N,Z_N}$ with
$Q_{\de_N}L_{\de_N} u_N =0$ in $\Omega\setminus B_{\de_N,Z_N}$ and
$\|u_N\|_\infty =1$ such that
\[
\|Q_{\de_N} L_{\de_N } u_N\|_{p} \le
\frac{1}{N}\frac{\de_N^{\frac{2}{p}}}{|\ln\de_N|^{\frac{(p-1)^2}{p}}}.
\]

First, we estimate $b^+_{i,h,N}$ and $b^-_{j,\bar h,N}$ in the
following formula:

\begin{equation}\label{3.6}
\begin{split}
Q_{\de_N} L_{\de_N } u_N= L_{\de_N } u_N&-\sum_{i=1}^m \sum_{h=1}^2
b^+_{i,h,N} \left(-\de_N^2\Delta\frac{\partial
P^+_{\de_N,Z_N,i}}{\partial z^+_{i,h}}\right)\\
&-\sum_{j=1}^n\sum_{\bar h=1}^2b^-_{j,\bar
h,N}\left(-\de_N^2\Delta\frac{\partial P^-_{\de_N,Z_N,j}}{\partial
z^-_{j,\bar h}}\right).
\end{split}
\end{equation}

For each fixed $k$, multiplying \eqref{3.6} by $
 \frac{\partial
P^+_{\de_N,Z_N,k}}{\partial z^+_{k,\hat h}}$, noting that

\[
 \int_\Om\bigl( Q_{\de_N} L_{\de_N } u_N\bigr)
 \frac{\partial
P^+_{\de_N,Z_N,k}}{\partial z^+_{k,\hat h}}=0,
\]
we obtain

\[
\begin{split}
& \int_\Om
  u_N\, L_{\de_N} \left(\frac{\partial
P^+_{\de_N,Z_N,k}}{\partial z^+_{k,\hat h}} \right)= \int_\Om\bigl(
  L_{\de_N} u_N\bigr)\, \frac{\partial
P^+_{\de_N,Z_N,k}}{\partial z^+_{k,\hat h}} \\
&=\sum_{i=1}^m
 \sum_{ h=1}^2 b^+_{i,h,N} \int_\Om\left(-\de_N^2\Delta\frac{\partial P^+_{\de_N,Z_N,i}}{\partial
 z^+_{i,h}}\right)
\frac{\partial P^+_{\de_N,Z_N,k}}{\partial z^+_{k,\hat
h}}\\
&\quad+\sum_{j=1}^n
 \sum_{\bar h=1}^2 b^-_{j,\bar h,N} \int_\Om\left(-\de_N^2\Delta\frac{\partial P^-_{\de_N,Z_N,j}}{\partial
 z^-_{j,\bar h}}\right)
\frac{\partial P^+_{\de_N,Z_N,k}}{\partial z^+_{k,\hat h}}.
\end{split}
\]
Using \eqref{2.9}, \eqref{2.91}  and
 Lemma~\ref{al1}, we obtain

\[
\begin{split}
&\int_\Om
  u_N\, L_{\de_N} \left(\frac{\partial
P^+_{\de_N,Z_N,k}}{\partial z^+_{k,\hat h}}\right)\\
&=\int_\Om\Bigg[-\de_N^2\Delta \left(\frac{\partial
P^+_{\de_N,Z_N,k}}{\partial z^+_{k,\hat h}}\right)-\sum_{i=1}^m p
\chi_{\Om_i^+}\left(
P^+_{\de_N,Z_N}-P_{\de_N,Z_N}^--\kappa_i^+-\frac{2\pi
q(x)}{|\ln\ep_N|}\right)_+^{p-1}\frac{\partial
P^+_{\de_N,Z_N,k}}{\partial z^+_{k,\hat h}}\\
&\quad-\sum_{j=1}^np\chi_{\Om_j^-}\left(
P^-_{\de_N,Z_N}-P^+_{\de_N,Z_N}-1+\frac{2\pi
q(x)}{\kappa|\ln\ep_N|}\right)_+^{p-1}\frac{\partial
P^+_{\de_N,Z_N,k}}{\partial z^+_{k,\hat h}}\Bigg]u_N\\
&=p\int_\Om\left(W_{\de_N,z^+_{k,N},a^+_{\de_N,k}}-a^+_{\de_N,k}\right)_+^{p-1}\left(\frac{\partial
W_{\de_N,z^+_{k,N},a^+_{\de_N,k}}}{\partial z^+_{k,\hat
h}}-\frac{\partial
a^+_{\de_N,k}}{\partial z^+_{k,\hat h}}\right)u_N\\
&\quad-p\sum_{\alpha=1}^m\int_{\Omega_\alpha^+}\left(W_{\de_N,z^+_{\alpha,N},a^+_{\de_N,\alpha}}-a^+_{\de_N,\alpha}
+O\left(\frac{s^+_{N,\alpha}}{|\ln\ep_N|}\right)\right)_+^{p-1}\frac{\partial
P^+_{\de_N,Z_N,k}}{\partial
z^+_{k,\hat h}}u_N\\
&\quad-p\sum_{\beta=1}^n\int_{\Omega_\beta^-}\left(W_{\de_N,z^-_{\beta,N},a^-_{\de_N,\beta}}-a^-_{\de_N,\beta}
+O\left(\frac{s^-_{N,\beta}}{|\ln\ep_N|}\right)\right)_+^{p-1}\frac{\partial
P^+_{\de_N,Z_N,k}}{\partial
z^+_{k,\hat h}}u_N \\
&=0\left(\frac{\ep_N^2}{|\ln\ep_N|^p}\right).
\end{split}
\]

Using \eqref{3.5} and  \eqref{3.51}, we find that
$$
b^+_{i,h,N}=0\left(\ep_N^2|\ln\ep_N|\right).
$$
Similarly,
$$
b^-_{i,h,N}=0\left(\ep_N^2|\ln\ep_N|\right).
$$
 Therefore,
\[
\begin{split}
&\sum_{i=1}^m\sum_{h=1}^2b^+_{i,h,N}\left(-\de_N^2\Delta\frac{\partial
P^+_{\de_N,Z_N,i}}{\partial z^+_{i,h}}\right)+\sum_{j=1}^n\sum_{\bar
h=1}^2b^-_{j,\bar h,N}\left(-\de_N^2\Delta\frac{\partial
P^-_{\de_N,Z_N,j}}{\partial z^-_{j,\bar h}}\right)\\
&=p\sum_{i=1}^m\sum_{h=1}^2b^+_{i,h,N}\left(W_{\de_N,z^+_{i,N},a^+_{\de_N,i}}-a^+_{\de_N,i}\right)_+^{p-1}
\left(\frac{\partial W_{\de_N,z^+_{i,N},a^+_{\de_N,i}}}{\partial
z^+_{i,h}}-\frac{\partial a^+_{\de_N,i}}{\partial z^+_{i,h}}\right)\\
&\quad +p\sum_{j=1}^n\sum_{\bar h=1}^2b^-_{j,\bar
h,N}\left(W_{\de_N,z^-_{j,N},a^-_{\de_N,j}}-a^-_{\de_N,j}\right)_+^{p-1}
\left(\frac{\partial W_{\de_N,z^-_{j,N},a^-_{\de_N,j}}}{\partial
z^-_{j,\bar h}}-\frac{\partial a^-_{\de_N,j}}{\partial z^-_{j,\bar h}}\right)\\
&=0\left(\sum_{i=1}^m\sum_{h=1}^2
\frac{\ep_N^{\frac{2}{p}-1}|b^+_{i,h,N}|}{|\ln\ep_N|^{p}}\right)
+0\left(\sum_{j=1}^n\sum_{\bar h=1}^2  \frac{\ep_N^{\frac{2}{p}-1}|b^-_{j,\bar h,N}|}{|\ln\ep_N|^{p}}\right)\\
&=0\left(\frac{\ep_N^{\frac{2}{p}+1}}{|\ln\ep_N|^{p-1}}\right)\quad
\text{in}~~L^p(\Om).
\end{split}
\]

Thus, we obtain

\[
L_{\de_N}u_N = Q_{\de_N}L_{\de_N}u_N
+O\left(\frac{\ep_N^{\frac{2}{p}+1}}{|\ln\ep_N|^{p-1}}\right)
 =O\left(\frac{1}{N}\frac{\de_N^{\frac{2}{p}}}{|\ln\de_N|^{\frac{(p-1)^2}{p}}}\right).
\]

For any fixed $i, j$, define

\[
\tilde u^+_{i,N} (y)= u_N(s^+_{N,i} y+z^+_{i,N}),\quad \tilde
u^-_{j,N} (y)= u_N(s^-_{N,j} y+z^-_{j,N}).
\]

Let
\[
\begin{split}
\tilde L_N^\pm u=& -\Delta u
-\sum_{k=1}^mp\frac{(s_{N,i}^\pm)^2}{\de_N^2}\chi_{\Om_k^+}\left(P^+_{\de_N,Z_N}(s^\pm_{N,i}
y+z^\pm_{i,N})-P^-_{\de_N,Z_N}(s^\pm_{N,i}
y+z^\pm_{i,N})-\kappa_k^+-\frac{2\pi
q}{|\ln\ep_N|}\right)_+^{p-1}u\\
&-\sum_{l=1}^np\frac{(s_{N,i}^\pm)^2}{\de_N^2}\chi_{\Om_l^-}\left(P^-_{\de_N,Z_N}(s^\pm_{N,i}
y+z^\pm_{i,N})-P^+_{\de_N,Z_N}(s^\pm_{N,i}
y+z^\pm_{i,N})-\kappa_l^-+\frac{2\pi q}{|\ln\ep_N|}\right)_+^{p-1}u.
\end{split}
\]

 Then

\[
(s_{N,i}^\pm)^{\frac{2}{p}}\times\frac{\de_N^2}{(s^\pm_{N,i})^2}\|\tilde
L_N^\pm \tilde u^\pm_{i,N}\|_p= \| L_{\de_N} u_N\|_p.
\]

Noting that
$$ \left(\frac{\de_N}{s^\pm_{N,i}}\right)^2=O\left(\frac{1}{|\ln\de_N|^{p-1}}\right),$$
we find that
$$
L_{\de_N}u_N=o\left(\frac{\de_N^{\frac{2}{p}}}{|\ln\de_N|^{\frac{(p-1)^2}{p}}}\right).
$$
 As a result,

\[
\tilde  L_{N}^\pm \tilde u_{i,N}^\pm
 =o(1),\quad \text{in}\; L^p(\Omega_N^\pm),
\]
where $\Omega_N^\pm=\bigl\{y: s^\pm_{N,i}
y+z^\pm_{i,N}\in\Omega\bigr\}$.

Since $\|\tilde u_{i,N}^\pm\|_\infty=1$, by the regularity theory of
elliptic equations,  we may assume that

\[
\tilde u_{i,N}^\pm\to u_i^\pm,  \quad \text{in}\;
C_{loc}^1(\mathbb{R}^2).
\]

It is easy to see that
\[
\begin{split}
&\sum_{k=1}^m\frac{(s^+_{N,i})^2}{\de_N^2}\chi_{\Om_k^+}\left(
P^+_{\de_N,Z_N}(s^+_{N,i}y+z^+_{i,N})-P^-_{\de_N,Z_N}(s_{N,i}^+y+z_{i,N}^+)-\kappa_k^+-\frac{2\pi q}{|\ln\ep_N|}\right)_+^{p-1}\\
&=\frac{(s_{N,i}^+)^2}{\de_N^2}\left(W_{\de_N,z^+_{i,N},a^+_{\de_N,i}}-a^+_{\de_N,i}
+O\left(\frac{s^+_{N,i}}{|\ln\ep_N|}\right)\right)_+^{p-1}+o(1)\\
&\rightarrow w_+^{p-1}.
\end{split}
 \]
 Similarly,
 \[
\begin{split}
&\sum_{l=1}^n\frac{(s^-_{N,j})^2}{\de_N^2}\chi_{\Om_l^-}\left(
P^-_{\de_N,Z_N}(s_{N,j}^-y+z_{j,N}^-)-P^+_{\de_N,Z_N}(s^-_{N,j}y+z^-_{j,N})-\kappa_l^-+\frac{2\pi q}{|\ln\ep_N|}\right)_+^{p-1}\\
&\rightarrow w_+^{p-1}.
\end{split}
 \]

  Then, by Lemma~\ref{al1}, we find that
  $u_i^\pm$ satisfies

\[
-\Delta u-pw_+^{p-1} u= 0.
\]
Now from the Proposition 3.1 in \cite{DY}, we have

\begin{equation}\label{3.7}
u_i^\pm= c_1^\pm \frac{\partial w}{\partial x_1}+ c_2^\pm
\frac{\partial w}{\partial x_2}.
\end{equation}

Since

\[
\int_\Om \Delta \bigl(\frac{\partial P^\pm_{\de_N,Z_N,i}}{\partial
z^\pm_{i, h}}\bigr) u_N =0,
\]
we find that

\[
\int_{\mathbb{R}^2}\phi_+^{p-1} \frac{\partial \phi}{\partial z_h}
u_i^\pm =0,
\]
which, together with \eqref{3.7}, gives $u_i^\pm \neq 0$. Thus,

\[
\tilde u_{i,N}^\pm \to 0,\quad \text{in}\; C^1(B_{L}(0)),
\]
for any $L>0$, which implies that $u_N=o(1)$ on $\partial
B_{Ls^\pm_{N,i}}(z^\pm_{i,N})$.

By assumption,

\[
Q_{\de_N} L_{\de_N} u_N = 0,\quad\text{in}\; \Om\setminus
B_{\de_N,Z_N}.
\]

On the other hand, by Lemma~\ref{al1}, for $i=1,\cdots,m$,
$j=1,\cdots,n$, we have

\[
\left( P^+_{\de_N,Z_N}-P_{\de_N,Z_N}^- -\kappa_i^+-\frac{2\pi
q(x)}{|\ln\ep_N|} \right)_+ =0, \quad x\in \Om_i^+\setminus
 B_{L s^+_{N,i}}(z^+_{i,N}),
\]
\[
\left( P^-_{\de_N,Z_N}-P_{\de_N,Z_N}^+ -\kappa_j^-+\frac{2\pi
q(x)}{|\ln\ep_N|} \right)_+ =0, \quad x\in \Om_j^-\setminus B_{L
s^-_{N,j}}(z^-_{j,N}).
\]
Thus, we find that

\[
-\Delta u_N  =0,\quad \text{in}~\Om\setminus  B_{\de_N,Z_N}.
\]
However, $u_N=0$ on $\partial\Om$ and $u_N=o(1)$ on $\partial
B_{\de_N,Z_N}$. So we have

\[
u_N=o(1).
\]
 This is a contradiction.

\end{proof}

From Lemma \ref{l31}, using Fredholm alternative, we can prove, as
in \cite{CLW},  the following result:

\begin{proposition}\label{p32}

$Q_\de L_\de$ is one to one and onto from $E_{\de,Z}$ to
$F_{\de,Z}$.

\end{proposition}

Now  consider the equation

\begin{equation}\label{3.10}
Q_\de L_\de \omega=
 Q_\de l_\de^+- Q_\de l_\de^-+ Q_\de R^+_\de(\omega)-Q_\de R^-_\de(\omega),
\end{equation}
where

\begin{equation}\label{3.11}
l_\de^+ = \sum_{i=1}^m\chi_{\Om_i^+}\left( P_{\de,Z}^+-P^-_{\de,Z}
-\kappa_i^+-\frac{2\pi q(x)}{|\ln\ep|} \right)_+^{p}-\sum_{i=1}^m
\left(W_{\de,z_i^+,a^+_{\de,i}}-a^+_{\de,i}\right)_+^{p},
\end{equation}
\begin{equation}\label{3.111}
l_\de^- = \sum_{j=1}^n\chi_{\Om_j^-}\left( P_{\de,Z}^--P^+_{\de,Z}
-\kappa_j^-+\frac{2\pi q(x)}{|\ln\ep|} \right)_+^{p}-\sum_{j=1}^n
\left(W_{\de,z_j^-,a^-_{\de,j}}-a^-_{\de,j}\right)_+^{p},
\end{equation}
and

\begin{equation}\label{3.12}
\begin{split}
R^+_\de(\omega)=&
\sum_{i=1}^m\chi_{\Om_i^+}\Bigg[\left(P^+_{\de,Z}-P_{\de,Z}^-+\omega
-\kappa_i^+-\frac{2\pi q(x)}{|\ln\ep|} \right)_+^{p} -
 \left( P^+_{\de,Z}-P_{\de,Z}^-
-\kappa_i^+-\frac{2\pi q(x)}{|\ln\ep|}\right)_+^{p} \\
&- p\left(P^+_{\de,Z}-P_{\de,Z}^- -\kappa_i^+-\frac{2\pi
q(x)}{|\ln\ep|}\right)_+^{p-1}\omega\Bigg],
\end{split}
\end{equation}
\begin{equation}\label{3.121}
\begin{split}
R^-_\de(\omega)=&
\sum_{j=1}^n\chi_{\Om_j^-}\Bigg[\left(P^-_{\de,Z}-P_{\de,Z}^+-\omega
-\kappa_j^-+\frac{2\pi q(x)}{|\ln\ep|} \right)_+^{p} -
 \left( P^-_{\de,Z}-P_{\de,Z}^+
-\kappa_j^-+\frac{2\pi q(x)}{|\ln\ep|}\right)_+^{p} \\
&+p\left(P^-_{\de,Z}-P_{\de,Z}^+ -\kappa_j^-+\frac{2\pi
q(x)}{|\ln\ep|}\right)_+^{p-1}\omega\Bigg].
\end{split}
\end{equation}

Using Proposition~\ref{p32}, we can rewrite \eqref{3.10} as

\begin{equation}\label{3.13}
\omega =G_\de\omega =: (Q_\de L_\de)^{-1} Q_\de \bigl(
  l_\de^+-l_\de^-+ R^+_\de(\omega)-R^-_\de(\omega)\bigr).
\end{equation}

The next Proposition enables us to reduce the problem of finding a
solution for \eqref{1} to a finite dimensional problem.

\begin{proposition}\label{p33}

There is an $\de_0>0$, such that for any $\de\in (0,\de_0]$ and  $Z$
 satisfying  \eqref{2.5},  \eqref{3.10} has a unique solution $\omega_\de\in
 E_{\de,Z}$, with

\[
\|\omega_\de\|_\infty =0\Bigl(\de|\ln\de|^{\frac{p-1}{2}}\Bigr).
\]
\end{proposition}

\begin{proof}

It follows from Lemma~\ref{al1} that if $L$ is large enough, $\de$
is small then

\[
 \left(P^+_{\de,Z}-P^-_{\de,Z}
-\kappa_i^+-\frac{2\pi q(x)}{|\ln\ep|} \right)_+=0, \quad x\in
\Om_i^+\setminus
 B_{L s^+_{\de,i}}(z^+_{i}),\,i=1,\cdots,m
\]
\[
 \left( P^-_{\de,Z}-P_{\de,Z}^+ -\kappa_j^-+\frac{2\pi
q(x)}{|\ln\ep|} \right)_+ =0, \quad x\in \Om_j^-\setminus B_{L
s^-_{\de,j}}(z^-_{j}),\,j=1,\cdots,n.
\]

Let

\[
M=  E_{\de,Z}\cap\Bigl\{ \|\omega\|_\infty\le
\de|\ln\de|^{\frac{p-1}{2}}\Big\}.
\]
Then $M$ is complete under $L^\infty$ norm  and $G_\de$ is a map
from $ E_{\de,Z}$ to $ E_{\de,Z}$. We will show that $G_\de $ is a
contraction map from $M$ to $M$.

Step~1.  $G_\de$ is a map from $M$ to $M$.

 For any $\omega\in M$, similar to Lemma~\ref{al1},
it is easy to  prove that for large $L>0$, $\de$ small

\begin{equation}\label{3.14}
\begin{split}
 \left(P^+_{\de,Z}-P^-_{\de,Z}+\omega
-\kappa_i^+-\frac{2\pi q(x)}{|\ln\ep|} \right)_+=0,  \quad x\in
\Om_i^+\setminus
 B_{L s^+_{\de,i}}(z^+_{i}),\\
 \left(P^-_{\de,Z}-P^+_{\de,Z}-\omega
-\kappa_j^-+\frac{2\pi q(x)}{|\ln\ep|} \right)_+=0, \quad x\in
\Om_j^-\setminus B_{L s^-_{\de,j}}(z^-_{j}).
\end{split}
\end{equation}
Note also that for any $u\in L^\infty(\Om)$,

\[
Q_\de u= u\quad \text{in}\; \Om\setminus B_{\de,Z}.
\]
Therefore, using Lemma~\ref{al1}, \eqref{3.11}--\eqref{3.121}, we
find that for any $\omega\in M$,
\[
\begin{split}
&Q_\de( l_\de^+-l_\de^-) + Q_\de( R_\de^+(\omega)-R_\de^-(\omega))\\
 =&l_\de^+-l_\de^-+R_\de^+(\omega)-R_\de^-(\omega)\\
 =&0, \quad \text{in}\; \Om\setminus  B_{Z,\de}.
 \end{split}
\]
So, we can apply Lemma~\ref{l31} to obtain

\[
\begin{split}
& \| (Q_\de L_\de)^{-1} \bigl(
 Q_\de( l_\de^+-l_\de^-) + Q_\de(
 R_\de^+(\omega)-R_\de^-(\omega))\bigr)\|_\infty\\
\le & \frac{C|\ln\de|^{\frac{(p-1)^2}{p}}}{\de^{\frac{2}{p}}}\|
Q_\de( l_\de^+-l_\de^-) + Q_\de(
R_\de^+(\omega)-R_\de^-(\omega))\|_p.
 \end{split}
\]

Thus, for any
 $\omega\in M$, we have

\begin{equation}\label{3.15}
\begin{array}{ll}
\| G_\de(\omega)\|_\infty= \| (Q_\de L_\de)^{-1} Q_\de \bigl(
 l_\de^+-l_\de^- +  R_\de^+(\omega)-R_\de^-(\omega)\bigr)\|_\infty&\\
 \qquad\quad\quad\quad\le  \frac{C|\ln\de|^{\frac{(p-1)^2}{p}}}{\de^{\frac{2}{p}}}\|Q_\de \bigl(  l_\de^+-l_\de^- +
R_\de^+(\omega)-R_\de^-(\omega)\bigr)\|_p.
\end{array}
\end{equation}

It follows from \eqref{3.4}--\eqref{3.51} that the constant
$b_{k,\hat h}^\pm$, corresponding to $u\in L^\infty(\Om)$, satisfies

\[
|b_{k,\hat h}^\pm| \le C |\ln\de|^{p+1}\Bigg(\sum_{i,\, h}
 \int_\Om \Bigl|\frac{\partial P^+_{\de,Z,i}}
{\partial z^+_{i, h}}\Bigr| |u|+\sum_{j,\, \bar h}
 \int_\Om \Bigl|\frac{\partial P^-_{\de,Z,j}}
{\partial z^-_{j, \bar h}}\Bigr| |u|\Bigg).
\]

 Since

\[
l_\de^+-l_\de^- +  R_\de^+(\omega)-R_\de^-(\omega)=0, \quad
\text{in}\; \Om\setminus B_{\de,Z},
\]
 we find that the constant
$b_{k,\hat h}^\pm$, corresponding to $l_\de^+-l_\de^- +
R_\de^+(\omega)-R_\de^-(\omega)$ satisfies

\[
\begin{split}
|b_{k,\hat h}^\pm|\le & C|\ln\de|^{p+1} \sum_{i,\, h}
\left(\sum_{\alpha=1}^m \int_{B_{Ls^+_{\de,\alpha} }(z_\alpha^+)}
 \Bigl|\frac{\partial P^+_{\de,Z,i}}
{\partial z^+_{i, h}}\Bigr||l_\de^+-l_\de^- +
R_\de^+(\omega)-R_\de^-(\omega)
|\right)\\
& +  C|\ln\de|^{p+1} \sum_{j,\,\bar h} \left(\sum_{\beta=1}^n
\int_{B_{Ls^-_{\de,\beta} }(z_\beta^-)}
 \Bigl|\frac{\partial P^-_{\de,Z,j}}
{\partial z^-_{j, \bar h}}\Bigr||l_\de^+-l_\de^- +
R_\de^+(\omega)-R_\de^-(\omega)
|\right)\\
\le & C\ep^{1-\frac{2}{p}}|\ln\ep|^{p}\|l_\de^+-l_\de^- +
R_\de^+(\omega)-R_\de^-(\omega) \|_p.
\end{split}
\]
As a result,

\[
\begin{split}
&\|Q_\de (l_\de^+-l_\de^- +  R_\de^+(\omega)-R_\de^-(\omega))\|_p \\
\le & \| l_\de^+-l_\de^- +
R_\de^+(\omega)-R_\de^-(\omega)\|_p+C \sum_{i,\,h} |b_{i,h}^+|
\left\| -\de^2\Delta\Bigl(\frac{\partial P^+_{\de,Z,i}}{\partial
z^+_{i,h}}\Bigr)
\right\|_p\\
&+ C \sum_{j,\,\bar h} |b_{j,\bar h}^-| \left\|
-\de^2\Delta\Bigl(\frac{\partial P^-_{\de,Z,j}}{\partial
z^-_{j,\bar h}}\Bigr) \right\|_p\\
\le & C\bigl( \| l^+_\de\|_p+ \| l^-_\de\|_p+ \|
R^+_\de(\omega)\|_p+\| R^-_\de(\omega) \|_p\bigr).
\end{split}
\]

On the other hand, from Lemma~\ref{al1} and \eqref{2.9}, we can
deduce

\[
\begin{split}
\|l_\de^+\|_p =&\left\|\sum_{i=1}^m\chi_{\Om_i^+}\left(
P_{\de,Z}^+-P^-_{\de,Z} -\kappa_i^+-\frac{2\pi q(x)}{|\ln\ep|}
\right)_+^{p}-\sum_{i=1}^m
\left(W_{\de,z_i^+,a^+_{\de,i}}-a^+_{\de,i}\right)_+^{p}\right\|_p\\
\le&\sum_{i=1}^m\frac{Cs^+_{\de,i}}{|\ln\ep|}\Big\|\bigl(W_{\de,z_i^+,a^+_{\de,i}}-a^+_{\de,i}\bigr)_+^{p-1}\Big\|_p\\
=&O\left(\frac{\de^{1+\frac{2}{p}}}{|\ln\de|^{\frac
{p-1}2+\frac{1}{p}}}\right).
\end{split}
\]

For the estimate of $\| R^+_\de(\omega) \|_p$, we have
\begin{equation}\label{3.16}
\begin{split}
\| R^+_\de(\omega)
\|_\infty=&\bigg\|\sum_{i=1}^n\chi_{\Om_i^+}\bigg[\Big(
P_{\de,Z}^+-P_{\de,Z}^-+\omega-\kappa_i^+-\frac{2\pi
q(x)}{|\ln\ep|}\Big)_+^{p}-\Big(P_{\de,Z}^+-P_{\de,Z}^--\kappa_i^+-\frac{2\pi q(x)}{|\ln\ep|}\Big)_+^{p}\\
&-p\Big(P_{\de,Z}^+-P_{\de,Z}^--\kappa_i^+-\frac{2\pi q(x)}{|\ln\ep|}\Big)_+^{p-1}\omega\bigg]\bigg\|_p\\
\le&C\|\omega\|_\infty^2\left\|\sum_{i=1}^n\chi_{\Om_i^+}\left(P_{\de,Z}^+-P_{\de,Z}^--\kappa_i^+-\frac{2\pi q(x)}{|\ln\ep|}\right)_+^{p-2}\right\|_p\\
=&O\left(\frac{\de^{\frac{2}{p}}\|\omega\|_\infty^2}{|\ln\de|^{p-3+\frac{1}{p}}}\right).
\end{split}
\end{equation}
Similarly, we have
\[
\|l_\de^-\|_p=O\left(\frac{\de^{1+\frac{2}{p}}}{|\ln\de|^{\frac
{p-1}2+\frac{1}{p}}}\right) ,\quad \| R^-_\de(\omega)
\|_p=O\left(\frac{\de^{\frac{2}{p}}\|\omega\|_\infty^2}{|\ln\de|^{p-3+\frac{1}{p}}}\right).
\]

Thus, we obtain

\begin{equation}\label{3.17}
\begin{split}
\| G_\de(\omega)\|_\infty \le &
\frac{C|\ln\de|^{\frac{(p-1)^2}{p}}}{\de^{\frac{2}{p}}}\Bigl(\|
l^+_\de \|_p+\| l^-_\de \|_p+\| R^+_\de(\omega)\|_p+\|
R^-_\de(\omega)\|_p
\Bigr)\\
\le  & C|\ln\de|^{\frac{(p-1)^2}{p}}\left(\frac{\de}{|\ln\de|^{\frac
{p-1}2+\frac{1}p}}+\frac{\|\omega\|_\infty^2}{|\ln\de|^{p-3+\frac{1}p}}\right)\\
\le&\de|\ln\de|^{\frac{p-1}2}
\end{split}
\end{equation}

Thus, $G_\de$ is a map from $M$ to $M$.

 Step~2.  $G_\de$ is a
contraction map.

In fact, for any $\omega_i\in M$, $i=1,2 $, we have

\[
G_\de \omega_1-G_\de \omega_2= (Q_\de L_\de)^{-1} Q_\de \bigl[
R^+_\de(\omega_1)-R^+_\de(\omega_2)-(R^-_\de(\omega_1)-R^-_\de(\omega_2))\bigr].
\]
Noting that

\[
R^+_\de(\omega_1)=R^+_\de(\omega_2)=0,\quad\text{in}\; \Om\setminus
\cup_{i=1}^m B_{Ls^+_{\de,i}}(z_i^+),
\]
and
\[
R^-_\de(\omega_1)=R^-_\de(\omega_2)=0,\quad\text{in}\; \Om\setminus
\cup_{j=1}^n B_{Ls^-_{\de,j}}(z_j^-),
\]
 we can deduce as in Step~1 that
\[
\begin{split}
\|G_\de \omega_1-G_\de \omega_2\|_\infty\le& \frac{C|\ln\de|^{\frac{(p-1)^2}{p}}}{\de^{\frac{2}{p}}}
(\|R^+_\de(\omega_1)-R^+_\de(\omega_2)\|_p+\|R^-_\de(\omega_1)-R^-_\de(\omega_2)\|_p)\\
\le&C|\ln\de|^{p-1}\left(\frac{\|\omega_1\|_\infty}{|\ln\de|^{p-2}}+\frac{\|\omega_2\|_\infty}{|\ln\de|^{p-2}}\right)\|\omega_1-\omega_2\|_\infty\\
\le&C\de|\ln\de|^{\frac {p+1}2}\|\omega_1-\omega_2\|_\infty \le\frac
12 \|\omega_1-\omega_2\|_\infty.
\end{split}
\]

Combining Step~1 and Step~2,  we have proved that $G_\de$ is a
contraction map from $M$ to $M$. By
  the
 contraction mapping theorem, there is an unique $\omega_\de\in M$, such
that $\omega_\de= G_\de\omega_\de$. Moreover, it follows from
\eqref{3.17} that

\[
\|\omega_\de\|_\infty\le \de|\ln\de|^{\frac {p-1}2}.
\]

\end{proof}

\section{Proof of The main results}

In this section, we will choose $Z$, such that $
P^+_{\de,Z}-P^-_{\de,Z}+\omega_\de$, where $\omega_\de$ is the map
obtained in Proposition~\ref{p33}, is a solution of \eqref{1}.

Define
\[
\begin{split}
I(u)=\frac{\de^2}{2}\int_\Omega
|Du|^2&-\sum_{i=1}^m\frac{1}{p+1}\int_\Omega\chi_{\Om_i^+}\left(u-\kappa_i^+-\frac{2\pi
q(x)}{|\ln\ep|}\right)_+^{p+1} \\
&-\sum_{j=1}^n\frac{1}{p+1}\int_\Omega\chi_{\Om_j^-}\left(\frac{2\pi
q(x)}{|\ln\ep|}-\kappa_j^--u\right)_+^{p+1}
\end{split}
\]
and
\begin{equation}\label{K}
 K(Z)= I\left(P_{\de,Z}^+-P^-_{\de,Z} +\omega_\de\right).
\end{equation}
It is well known that  if $Z$ is a critical point of $K(Z)$, then
$P_{\de,Z}^+-P^-_{\de,Z} +\omega_\de$ is a solution of \eqref{1}. In
the following,  we will prove that $K(Z)$ has a critical point.

 \begin{lemma}\label{l42}

 We have

 \[
 K(Z)= I\left(P^+_{\de,Z}-P^-_{\de,Z}
\right)+O\left(
 \frac{\ep^3}{|\ln\ep|^p}\right).
\]

 \end{lemma}

\begin{proof}

Recall that
$
P^+_{\de,Z}=\sum_{i=1}^m P^+_{\de,Z,i},\quad
P^-_{\de,Z}=\sum_{j=1}^n P^-_{\de,Z,j}.
$
We have
\[
\begin{split}
K(Z)&= I\bigl(P^+_{\de,Z}-P^-_{\de,Z}\bigr) +\de^2\int_\Om
D\bigl(P^+_{\de,Z}-P^-_{\de,Z}\bigr)D\omega_\de+\frac{\de^2}{2}\int_\Omega |D\omega_\de|^2\\
&\quad\quad-\sum_{i=1}^m\frac1{p+1} \int_\Om \chi_{\Om_i^+}\Biggl[
\biggl( P^+_{\de,Z}-P^-_{\de,Z} +\omega_\de-\kappa_i^+-\frac{2\pi
q(x)}{|\ln\ep|} \biggr)_+^{p+1}\\
&\quad\quad\quad- \biggl(
P_{\de,Z}^+-P_{\de,Z}^--\kappa_i^+-\frac{2\pi q(x)}{|\ln\ep|}
\biggr)_+^{p+1}\Biggr]\\
&\quad\quad-\sum_{j=1}^n\frac1{p+1} \int_\Om \chi_{\Om_j^-}\Biggl[
\biggl( P^-_{\de,Z}-P^+_{\de,Z} -\omega_\de-\kappa_j^-+\frac{2\pi
q(x)}{|\ln\ep|} \biggr)_+^{p+1}\\
&\quad\quad\quad- \biggl(
P_{\de,Z}^--P_{\de,Z}^+-\kappa_j^-+\frac{2\pi q(x)}{|\ln\ep|}
\biggr)_+^{p+1}\Biggr].
\end{split}
\]

Using Proposition~\ref{p33} and \eqref{3.14}, we find

\[
\begin{split}
&  \int_{\Om_i^+} \Biggl[ \biggl( P^+_{\de,Z}-P^-_{\de,Z}
+\omega_\de-\kappa_i^+-\frac{2\pi q(x)}{|\ln\ep|} \biggr)_+^{p+1}-
\biggl( P_{\de,Z}^+-P_{\de,Z}^--\kappa_i^+-\frac{2\pi
q(x)}{|\ln\ep|}
\biggr)_+^{p+1}\Biggr]\\
= &\int_{ B_{Ls^+_{\de,i}}(z^+_i)} \Biggl[ \biggl(
P^+_{\de,Z}-P^-_{\de,Z} +\omega_\de-\kappa_i^+-\frac{2\pi
q(x)}{|\ln\ep|} \biggr)_+^{p+1}- \biggl(
P_{\de,Z}^+-P_{\de,Z}^--\kappa_i^+-\frac{2\pi q(x)}{|\ln\ep|}
\biggr)_+^{p+1}\Biggr]\\
=&O\left(\frac{(s^+_{\de,i})^2\|\omega_\de\|_\infty}{|\ln\ep|^{p}}\right)\\
= &O\Bigl(
 \frac{\ep^3}{|\ln\ep|^{p}}\Bigr).
\end{split}
\]

On the other hand,

\[
\begin{split}
&\de^2 \int_\Om DP^+_{\de,Z}D\omega_\de
 =
 \sum_{i=1}^m\int_\Om
\left(W_{\de,z^+_i,a^+_{\de,i}}-a^+_{\de,i}\right)_+^{p} \omega_\de\\
=&\sum_{i=1}^m
 \int_{ B_{s^+_{\de,k}}(z^+_k)}
(W_{\de,z_i^+,a^+_{\de,i}}-a^+_{\de,i})_+^{p} \omega_\de \\
=& O\Bigl(
 \frac{\ep^3}{|\ln\ep|^{p}}\Bigr).
\end{split}
\]
Next, we estimate $\de^2\int_\Om |D\omega_\de|^2$. Note that
\[
\begin{split}
-\de^2\Delta\omega_\de=&\sum_{i=1}^m\chi_{\Om_i^+}\left(P^+_{\de,Z}-P^-_{\de,Z}+\omega_\de-\kappa_i^+-\frac{2\pi
q(x)}{|\ln\ep|}\right)_+^{p}
-\sum_{i=1}^m\left(W_{\de,z_i^+,a^+_{\de,i}}-a^+_{\de,i}\right)_+^{p}\\
&-\sum_{j=1}^n\chi_{\Om_j^-}\left(P^-_{\de,Z}-P^+_{\de,Z}-\omega_\de-\kappa_j^-+\frac{2\pi
q(x)}{|\ln\ep|}\right)_+^{p}+\sum_{j=1}^n\left(W_{\de,z_j^-,a^-_{\de,j}}-a^-_{\de,j}\right)_+^{p}\\
&+\sum_{i=1}^m\sum_{ h=1}^2b_{i,
h}^+\left(-\de^2\Delta\frac{\partial P^+_{\de,Z,i}}{\partial
z^+_{i,h}}\right)+\sum_{j=1}^n\sum_{\bar h=1}^2b^-_{j, \bar
h}\left(-\de^2\Delta\frac{\partial P^-_{\de,Z,j}}{\partial
z^-_{j,\bar h}}\right).
\end{split}
\]
Hence, by \eqref{2.9}--\eqref{2.91}, we have
\[
\begin{split}
\de^2\int_\Om|D\omega_\de|^2
=&\sum_{i=1}^m\int_{\Om_i^+}\Biggl[\left(P^+_{\de,Z}-P^-_{\de,Z}+\omega_\de-\kappa_i^+-\frac{2\pi
q(x)}{|\ln\ep|}\right)_+^{p}
-\left(W_{\de,z_i^+,a_{\de,i}^+}-a^+_{\de,i}\right)_+^{p}\Biggr]\omega_\de\\
&-\sum_{j=1}^n\int_{\Om_j^-}\Biggl[\left(P^-_{\de,Z}-P^+_{\de,Z}-\omega_\de-\kappa_j^-+\frac{2\pi
q(x)}{|\ln\ep|}\right)_+^{p}
-\left(W_{\de,z_j^-,a_{\de,j}^-}-a^-_{\de,j}\right)_+^{p}\Biggr]\omega_\de\\
&+\sum_{i=1}^m\sum_{ h=1}^2b^+_{i,
h}\int_\Om\left(-\de^2\Delta\frac{\partial P^+_{\de,Z,i}}{\partial
z^+_{i, h}}\right)\omega_\de+\sum_{j=1}^n\sum_{\bar
h=1}^2b^-_{j,\bar h}\int_\Om\left(-\de^2\Delta\frac{\partial
P^-_{\de,Z,j}}{\partial
z^-_{j,\bar h}}\right)\omega_\de\\
=&p\sum_{i=1}^m\int_{\Om_i^+}\left(W_{\de,z_i^+,a^+_{\de,i}}-a^+_{\de,i}\right)_+^{p-1}\left(\frac{s^+_{\de,i}}{|\ln\ep|}+\omega_\de\right)\omega_\de
+0\left(\sum_{i=1}^m\sum_{ h=1}^2\frac{\ep|b^+_{i,
h}|\|\omega_\de\|_\infty}{|\ln\ep|^{p}}\right)\\
&-p\sum_{j=1}^n\int_{\Om_j^-}\left(W_{\de,z^-_j,a^-_{\de,j}}-a^-_{\de,j}\right)_+^{p-1}\left(\frac{s^-_{\de,j}}{|\ln\ep|}+\omega_\de\right)\omega_\de
+0\left(\sum_{j=1}^n\sum_{\bar h=1}^2\frac{\ep|b^-_{j,\bar
h}|\|\omega_\de\|_\infty}{|\ln\ep|^{p}}\right)\\
 =&O\left(\frac{\ep^4}{|\ln\ep|^{p-1}}\right).
\end{split}
\]

Other terms can be estimated as above. So our assertion follows.

\end{proof}

\begin{lemma}\label{l43}

 We have

\[
\frac{\partial K(Z)}{\partial z^+_{i,h}}=\frac{\partial }{\partial
z^+_{i,h}}I\left(P^+_{\de,Z}-P^-_{\de,Z} \right)+O\Bigl(
\frac{\ep^3}{|\ln\ep|^{p-1}}\Bigr),\,\,i=1,\cdots,m,
\]
\[
\frac{\partial K(Z)}{\partial z^-_{j,\bar h}}=\frac{\partial
}{\partial z^-_{j,\bar h}}I\left(P^+_{\de,Z}-P^-_{\de,Z}
\right)+O\Bigl( \frac{\ep^3}{|\ln\ep|^{p-1}}\Bigr),\,\,j=1,\cdots,n.
\]
\end{lemma}

\begin{proof} We only give the proof of the first estimate.

First, we have

\begin{equation*}\label{4.4}
\begin{split}
&\frac{\partial K(Z)}{\partial z^+_{i,h}}=\left\langle
I^\prime\Bigl( P^+_{\de,Z}-P_{\de,Z}^-
+\omega_\de\Bigr),\frac{\partial P^+_{\de,Z}}{\partial
z^+_{i,h}}-\frac{\partial P^-_{\de,Z}}{\partial z^+_{i,h}}+
\frac{\partial
\omega_\de}{\partial z^+_{i,h}}\right\rangle\\
=&\frac{\partial }{\partial z^+_{i,h}}I\Bigl(P^+_{\de,Z}-P_{\de,Z}^-
\Bigr) +\left\langle I^\prime\big(
P^+_{\de,Z}-P_{\de,Z}^-+\omega_\de\big),
\frac{\partial\omega_\de}{\partial z^+_{i,h}}\right\rangle
\\
&-\sum_{k=1}^m \int_{\Omega_k^+}\Biggl[
\biggl(P^+_{\de,Z}-P_{\de,Z}^-+\omega_\de-\kappa_k^+-\frac{2\pi
q(x)}{|\ln\ep|}\biggr)_+^{p}-\biggl(P^+_{\de,Z}-P^-_{\de,Z}-\kappa_k^+-\frac{2\pi
q(x)}{|\ln\ep|}\biggr)_+^{p}\Biggr]\\
&\qquad \times\left(\frac{\partial P^+_{\de,Z}}{\partial
z^+_{i,h}}-\frac{\partial P^-_{\de,Z}}{\partial z^+_{i,h}}\right)\\
&- \sum_{l=1}^n\int_{\Omega_l^-}\Biggl[
\biggl(P^-_{\de,Z}-P_{\de,Z}^+-\omega_\de-\kappa_l^-+\frac{2\pi
q(x)}{|\ln\ep|}\biggr)_+^{p}-\biggl(P^-_{\de,Z}-P^+_{\de,Z}-\kappa_l^-+\frac{2\pi
q(x)}{|\ln\ep|}\biggr)_+^{p}\Biggr]\\
&\qquad \times\left(\frac{\partial P^-_{\de,Z}}{\partial
z^+_{i,h}}-\frac{\partial P^+_{\de,Z}}{\partial z^+_{i,h}}\right).
\end{split}
\end{equation*}

Since $\omega_\de\in E_{\de,Z}$, we have
\[
\int_\Om\left(W_{\de,z_k^\pm,a_{\de,k}^\pm}-a_{\de,k}^\pm\right)_+^{p-1}\left(\frac{\partial
W_{\de,z_k^\pm,a^\pm_{\de,k}}}{\partial z^\pm_{k,h}}-\frac{\partial
a_{\de,k}^\pm}{\partial z_{k,h}^\pm}\right)\omega_\de=0.
\]
Differentiating the above relation with respect to $z^+_{i, h}$, we
can deduce

\begin{equation*}
\begin{split}
&\Bigg\langle  I^\prime\bigl(P^+_{\de,Z}-P_{\de,Z}^-
+\omega_\de\bigr), \frac{\partial \omega_\de}{\partial z^+_{i, h}}
\Bigg\rangle\\
=&\sum_{\alpha=1}^m\sum_{\hat h=1}^2 b^+_{\alpha,\hat h}
\int_\Om\left(-\de^2\Delta\frac{\partial
P^+_{\de,Z,\alpha}}{\partial z^+_{\alpha,\hat h}}\right)
\frac{\partial \omega_\de}{\partial z^+_{i,
h}}+\sum_{\beta=1}^n\sum_{\tilde h=1}^2 b^-_{\beta,\tilde h}
\int_\Om\left(-\de^2\Delta\frac{\partial P^-_{\de,Z,\beta}}{\partial
z^-_{\beta,\tilde h}}\right) \frac{\partial \omega_\de}{\partial z^+_{i, h}}\\
=&\sum_{\alpha=1}^m\sum_{\hat h=1}^2 p b_{\alpha,\hat h}^+
\int_\Om\left(W_{\de,z_\alpha^+,a^+_{\de,\alpha}}-a_{\de,\alpha}^+\right)_+^{p-1}\left(\frac{\partial
W_{\de,z_\alpha^+,a^+_{\de,\alpha}}}{\partial z^+_{\alpha,\hat
h}}-\frac{\partial a^+_{\de,\alpha}}{\partial z^+_{\alpha,\hat
h}}\right) \frac{\partial
\omega_\de}{\partial z^+_{i, h}}\\
&+\sum_{\beta=1}^n\sum_{\tilde h=1}^2 p b_{\beta,\tilde h}^-
\int_\Om\left(W_{\de,z_\beta^-,a^-_{\de,\beta}}-a_{\de,\beta}^-\right)_+^{p-1}\left(\frac{\partial
W_{\de,z_\beta^-,a^-_{\de,\beta}}}{\partial z^-_{\beta,\tilde
h}}-\frac{\partial a^-_{\de,\beta}}{\partial z^-_{\beta,\tilde
h}}\right) \frac{\partial
\omega_\de}{\partial z^+_{i, h}}\\
 =&O\left(\sum_{\alpha=1}^m\sum_{\hat h=1}^2
 \frac{\ep|b^+_{\alpha,\hat
h}|}{|\ln\ep|^{p}}+\sum_{\beta=1}^n\sum_{\tilde h=1}^2
\frac{\ep|b^-_{\beta,\tilde
h}|}{|\ln\ep|^{p}}\right)=O\left(\frac{\ep^3}{|\ln\ep|^{p-1}}\right).
\end{split}
\end{equation*}
On the other hand, using \eqref{3.16} (for the definition of
$R^+_\de(\omega)$, see \eqref{3.12}), we obtain

\[
\begin{split}
&\sum_{k=1}^m\int_{\Omega_k^+}\Biggl[
\biggl(P^+_{\de,Z}-P_{\de,Z}^-+\omega_\de-\kappa_k^+-\frac{2\pi
q(x)}{|\ln\ep|}\biggr)_+^{p}-\biggl(P^+_{\de,Z}-P^-_{\de,Z}-\kappa_k^+-\frac{2\pi
q(x)}{|\ln\ep|}\biggr)_+^{p}\Biggr]\frac{\partial
P^+_{\de,Z,i}}{\partial z^+_{i,h}}\\
=&\sum_{k=1}^m \int_{\Omega_k^+}\Biggl[
\biggl(P^+_{\de,Z}-P_{\de,Z}^-+\omega_\de-\kappa_k^+-\frac{2\pi
q(x)}{|\ln\ep|}\biggr)_+^{p}
-\biggl(P^+_{\de,Z}-P_{\de,Z}^--\kappa_k^+-\frac{2\pi q(x)}{|\ln\ep|}\biggr)_+^{p}\\
&-p\biggl(P^+_{\de,Z}-P^-_{\de,Z}-\kappa_k^+-\frac{2\pi
q(x)}{|\ln\ep|}\biggr)_+^{p-1}\omega_\de\Biggr]\frac{\partial
P^+_{\de,Z,i}}{\partial z^+_{i,h}}\\
& +
\sum_{k=1}^mp\int_{\Omega_k^+}\Biggl[\biggl(P^+_{\de,Z}-P^-_{\de,Z}-\kappa_k^+-\frac{2\pi
q(x)}{|\ln\ep|}\biggr)_+^{p-1}-\bigl(W_{\de,z^+_k,a^+_{\de,k}}-a^+_{\de,k}\bigr)_+^{p-1}\Biggr]\frac{\partial
P^+_{\de,Z,i}}{\partial z^+_{i,h}}\omega_\de\\
&+O\left(\frac{(s^+_{\de,k})^2\|\omega_\de\|_\infty}{|\ln\ep|^{p}}\right)\\
=&\int_\Omega R^+_\de(\omega_\de)\frac{\partial
P^+_{\de,Z,i}}{\partial z^+_{i,h}} +\sum_{k=1}^m
p\int_{\Omega_k^+}\Biggl[\biggl(P^+_{\de,Z}-P_{\de,Z}^--\kappa_k^+-\frac{2\pi
q(x)}{|\ln\ep|}\biggr)_+^{p-1}\\
&-\bigl(W_{\de,z_k^+,a^+_{\de,k}}-a^+_{\de,k}\bigr)_+^{p-1}\Biggr]\frac{\partial
P^+_{\de,Z,i}}{\partial
z^+_{i,h}}\omega_\de+O\left(\frac{\ep^3}{|\ln\ep|^{p}}\right)\\
 =& O\left(\frac{\ep^3}{|\ln\ep|^{p-1}}\right).
\end{split}
\]
In addition, we have
\[
\begin{split}
&\int_{\Omega_l^+}\Biggl[
\biggl(P^+_{\de,Z}-P_{\de,Z}^-+\omega_\de-\kappa_l^+-\frac{2\pi
q(x)}{|\ln\ep|}\biggr)_+^{p}-\biggl(P^+_{\de,Z}-P^-_{\de,Z}-\kappa_l^+-\frac{2\pi
q(x)}{|\ln\ep|}\biggr)_+^{p}\Biggr]\frac{\partial
P^-_{\de,Z,i}}{\partial z^-_{i,h}}\\
&=p\int_{\Om_l^+}\biggl(P_{\de,Z}^+-P^-_{\de,Z}-\kappa_l^+-\frac{2\pi
q(x)}{|\ln\ep|}\biggr)_+^{p-1}\frac{\partial
P^-_{\de,Z,i}}{\partial z^-_{i,h}}\omega_\de\\
&=O\left(\frac{\ep^3}{|\ln\ep|^p}\right).
\end{split}
\]
Other teams can be estimated as above.  Thus, the estimate follows.
\end{proof}

\begin{proof}[Proof of Theorem~\ref{th3}]
Recall that $Z=(Z_m^+,Z_n^-).$ Set
 \[
\begin{split}
\Phi(Z_m^+,Z_n^-)=&\sum_{i=1}^m4\pi^2\kappa_i^+q(z_i^+)-\sum_{j=1}^n4\pi^2\kappa_j^-q(z_j^-)
+\sum_{i=1}^m\pi(\kappa_i^+)^2 g(z_i^+,z_i^+)\\
&+\sum_{j=1}^n\pi (\kappa_j^-)^2g(z_j^-,z_j^-)-\sum_{i\neq
k}\pi\kappa_i^+\kappa_k^+\bar
G(z_i^+,z_k^+)-\sum_{j\neq l}\pi \kappa_j^-\kappa_l^-\bar G(z_l^-,z_j^-)\\
&+\sum_{i=1}^m\sum_{j=1}^n2\pi\kappa_i^+\kappa_j^-\bar
G(z_i^+,z_j^-).
\end{split}
  \]

 Note that the Kirchhoff--Routh function associated to the
vortex dynamics now is
\[
\begin{split}
\mathcal{W}(Z_m^+,Z_n^-)
 =&\frac{1}{2}\sum_{i,k=1,i\neq
k}^m\kappa_i^+\kappa
_k^+G(z_i^+,z_k^+)+\frac{1}{2}\sum^{n}_{j,l=1,j\neq
l}\kappa^-_j\kappa^-_l
G(z_j^-,z_l^-)\\
 &+\frac 12\sum_{i=1}^m(\kappa_i^+)^2H(z_i^+,z_i^+)+\frac
12\sum_{j=1}^n(\kappa_j^-)^2H(z_j^-,z_j^-)\\
&-\sum_{i=1}^m\sum_{j=1}^n\kappa_i^+\kappa_j^-G(z_i^+,z_j^-)
+\sum^{m}_{i=1}\kappa_i^+\psi_0(z_i^+)-\sum^{n}_{j=1}\kappa_j^-\psi_0(z_j^-).
\end{split}\]
 Recall that  $h(z_i,z_j)=-H(z_i,z_j)$, it is easy to check
that
\[
\Phi(Z_m^+,Z_n^-)=-4\pi^2\mathcal {W}(Z_m^+,Z_n^-)+\pi\ln
R\left(\sum_{i=1}^m(\kappa_i^+)^2+\sum_{j=1}^n(\kappa_j^-)^2\right).
\]
Hence, $\Phi(Z_m^+,Z_N^-)$ and $\mathcal{W}(Z_m^+,Z_N^-)$ possess
the same critical points.

By Lemma \ref{l42}, \ref{l43} and Proposition \ref{ap1}, \ref{ap2},
we have
 \[
 K(Z)= \frac{ C\de^2
 }{\ln\frac{R}{\ep}
 }
 +\frac{\pi(p-1)\de^2}{4(\ln\frac{R}{\ep})^2}\left(\sum_{i=1}^m(\kappa_i^+)^2+\sum_{j=1}^n(\kappa_j^-)^2\right)
 +\frac{\de^2}{|\ln\ep|^2}\Phi(Z)
+ 0\left(\frac{\de^2\ln|\ln\ep|}{|\ln\ep|^3} \right)
\]
and
\[
\frac{\partial K(Z)}{\partial
z^\pm_{i,h}}=\frac{\de^2}{|\ln\ep|^2}\frac{\partial
\Phi(Z)}{\partial
z^\pm_{i,h}}+O\left(\frac{\de^2\ln|\ln\ep|}{|\ln\ep|^3} \right).
\]
Thus, the existence of a $C^1$-stable critical point of
Kirchhoff-Routh function $\mathcal{W}(Z)$ implies that  $K(Z)$ has a
critical point.

Thus we get a solution $w_\de$ for \eqref{1}. Let
$u_\ep=\frac{|\ln\ep|}{2\pi}w_\de,\,\de=\ep\left(\frac{|\ln\ep|}{2\pi}\right)^{\frac{1-p}{2}}$,
it is not difficult to check that  $u_\ep$ has all the properties
listed in Theorem \ref{th3} and thus the proof of Theorem \ref{th3}
is complete.
\end{proof}

Now we are in the position to prove Theorem~\ref{nth1}.
\begin{proof}[Proof of Theorem~\ref{nth1}]

By Theorem \ref{th3}, we obtain that $u_\ep$ is a solution to
\eqref{0}.

Set
 \[
 \begin{split}
&\mathbf{v}_\ep=(\nabla
(u_\ep-q))^\bot,\qquad\omega_\ep=\nabla\times\mathbf{v}_\ep,\\
P_\ep=&\sum_{i=1}^m\frac{1}{p+1}\chi_{\Om_i^+}\left(u_\ep-q-\frac{\kappa_i^+|\ln\ep|}{2\pi}\right)_+^{p+1}\\
&+\sum_{j=1}^n\frac{1}{p+1}\chi_{\Om_j^-}\left(q-\frac{\kappa_j^-|\ln\ep|}{2\pi}-u_\ep\right)_+^{p+1}-\frac{1}{2}|\nabla
(u_\ep-q)|^2.
\end{split}
\]
Then  $(\mathbf{v}_\ep, P_\ep)$ forms a stationary solution for
problem \eqref{1.2}.

We now just need to verify as $\ep\rightarrow0$
\[ \int_\Omega
\omega_\ep\rightarrow \sum_{j=1}^m\kappa_j^+-\sum_{j=1}^n\kappa_j^-.
\]

By direct calculations, we find that
\[
\begin{split}
\int_\Omega\omega_\ep&=\sum_{i=1}^m\frac{1}{\ep^2}\int_\Omega\chi_{\Om_i^+}
\left(u_\ep-q-\frac{\kappa_i^+|\ln\ep|}{2\pi}\right)_+^p-\sum_{j=1}^n\frac{1}{\ep^2}\int_\Omega\chi_{\Om_j^-}
\left(q-\frac{\kappa_j^-|\ln\ep|}{2\pi}-u_\ep\right)_+^p\\
&=\sum_{i=1}^m\frac{|\ln\ep|^p}{(2\pi)^p \ep^2}\int_{\Omega_i^+}
\left(w_\de-\kappa_i^+-\frac{2\pi
q}{|\ln\ep|}\right)^p_+-\sum_{j=1}^n\frac{|\ln\ep|^p}{(2\pi)^p
\ep^2}\int_{\Omega_j^-}
\left(\frac{2\pi q}{|\ln\ep|}-\kappa_j^--w_\de\right)^p_+\\
&=\frac{|\ln\ep|^p}{(2\pi)^p\ep^2}\sum_{i=1}^m\int_{B_{Ls^+_{\de,i}(z^+_i)}}
\left(W_{\de,z^+_i,a^+_{\de,i}}-a^+_{\de,i}+O\Big(\frac{s^+_{\de,i}}{|\ln\ep|}\Big)\right)^p_+\\
&\quad-\frac{|\ln\ep|^p}{(2\pi)^p\ep^2}\sum_{j=1}^n\int_{B_{Ls^-_{\de,j}(z^-_j)}}
\left(W_{\de,z^-_j,a^-_{\de,j}}-a^-_{\de,j}+O\Big(\frac{s^-_{\de,j}}{|\ln\ep|}\Big)\right)^p_+\\
&=\sum_{i=1}^m\frac{
(s^+_{\de,i})^2|\ln\ep|^p}{(2\pi)^p\ep^2}\left(\frac{\de}{s^+_{\de,i}}\right)^{\frac{2p}{p-1}}\int_{B_1(0)}\phi^p\\
&\quad-\sum_{j=1}^n\frac{ (s^-_{\de,j})^2|\ln\ep|^p}{(2\pi)^p\ep^2}\left(\frac{\de}{s^-_{\de,j}}\right)^{\frac{2p}{p-1}}\int_{B_1(0)}\phi^p+o(1)\\
&=\sum_{i=1}^m\frac{ a^+_{\de,i}
|\ln\ep|}{\ln\frac{R}{s^+_{\de,i}}}-\sum_{j=1}^n\frac{ a^-_{\de,j}
|\ln\ep|}{\ln\frac{R}{s^-_{\de,j}}}+o(1)\\
&\rightarrow \sum_{j=1}^m\kappa_j^+-\sum_{j=1}^n\kappa_j^-, \quad
\text{as}~~\ep\rightarrow0.
\end{split}
\]
Therefore, the result follows.
\end{proof}

\begin{remark}
To regularize pairs of  vortices  with equi-strength $\kappa$, we do
not need $\chi_{\Om_i^+}$ and $\chi_{\Om_j^-}$, that is, we only
 need to consider the following problem
\[
 \begin{cases}
-\ep^2 \Delta u=(u-q-\frac{\kappa}{2\pi}\ln\frac{1}{\ep})_+^p-(q-\frac{\kappa}{2\pi}\ln\frac{1}{\ep}-u)_+^p, \quad & x\in\Omega, \\
u=0, \quad & x\in\partial\Omega.
\end{cases}
\]

\end{remark}

{\bf Acknowledgements:} D. Cao and Z. Liu were supported by the
National Center for Mathematics and Interdisciplinary Sciences, CAS.
D. Cao and J. Wei were also supported by CAS Croucher Joint
Laboratories Funding Scheme.

\appendix

\section{Energy expansion}

In this section we will give precise expansions of $I\left(
P^+_{\de,Z}-P^-_{\de,Z}\right)$ and $\frac{\partial }{\partial
z^\pm_{i,h}} I\left( P_{\de,Z}^+-P_{\de,Z}^-\right)$, which have
been used in section 4.



We always assume that  $z_i^+, z_j^-\in\Om$ satisfies
\begin{equation*}
\begin{split}
&d(z_i^+,\partial\Om)\ge \varrho,~d(z_j^-,\partial\Om)\ge
\varrho,\quad |z_i^+-z_k^+|\ge \varrho^{\bar L},\quad i, k=1,\cdots,m,\; i\ne k\\
&~~|z_j^--z_l^-|\ge \varrho^{\bar L},\quad|z_i^+-z_j^-|\ge
\varrho^{\bar L},\quad j, l=1,\cdots,n,\; j\ne l,
\end{split}
\end{equation*}
where $\varrho>0$ is a fixed small constant and $\bar L>0$ is a
fixed large constant.

\begin{lemma}\label{al1}
For $x\in\Om_i^+,\,i=1,2,\cdots,m$ and
$x\in\Om_j^-,\,j=1,2,\cdots,m$, we have
\[
 P^+_{\de,Z}(x)-P^-_{\de,Z}(x)>\kappa_i^++\frac{2\pi q(x)}{|\ln\ep|},\quad
x\in B_{s^+_{\de,i}(1-Ts^+_{\de,i})}(z^+_i),
\]
\[
 P^-_{\de,Z}(x)-P^+_{\de,Z}(x)>\kappa_j^--\frac{2\pi q(x)}{|\ln\ep|},\quad
x\in B_{s^-_{\de,j}(1-Ts^-_{\de,j})}(z_j^-),
\]
where $T>0$ is a large constant; while

\[
 P_{\de,Z}^+(x)-P_{\de,Z}^-(x)<\kappa_i^++\frac{2\pi q(x)}{|\ln\ep|},\quad
x\in \Om_i^+\setminus
B_{s^+_{\de,i}(1+(s^+_{\de,i})^{\sigma})}(z^+_i),
\]
\[
 P_{\de,Z}^-(x)-P_{\de,Z}^+(x)<\kappa_j^--\frac{2\pi q(x)}{|\ln\ep|},\quad
x\in\Om_j^-\setminus
B_{s^-_{\de,j}(1+(s^-_{\de,j})^{\sigma})}(z^-_j),
\]
where $\sigma>0$ is a small constant.
\end{lemma}

\begin{proof}
The proof is exactly same as Lemma A.1 in \cite{CLW}. For reader's
convenience, we give the proof for $P_{\de,Z}^+-P_{\de,Z}^-$ here.

 Suppose that $x\in
B_{s^+_{\de,i}(1-Ts^+_{\de,i})}(z_i^+)$. It follows from \eqref{2.9}
and $\phi_1'(s)<0$ that
\[
\begin{split}
&P^+_{\de,Z}(x)-P_{\de,Z}^-(x)-\kappa_i^+-\frac{2\pi q(x)}{|\ln\ep|}=  W_{\de,z^+_i,a^+_{\de,i}}(x)-a^+_{\de,i}+O\left(\frac{s^+_{\de,i}}{|\ln\ep|}\right)\\
=&
\frac{a^+_{\de,i}}{|\phi'(1)||\ln\frac{R}{s^+_{\de,i}}|}\phi\Bigl(\frac{|x-z^+_i|}{s^+_{\de,i}}\Bigr)+O\Bigl(\frac{\ep}{|\ln\ep|}\Bigr)>0,
\end{split}
\]
if $T>0$ is large.   On the other hand, if $x\in\Om_i^+\setminus
 B_{(s^+_{\de,i})^{\tilde\sigma}}(z_i^+)$, where
$\tilde\sigma>\sigma>0$ is a fixed small constant,  then
\[
\begin{split}
& P^+_{\de,Z}(x)-P^-_{\de,Z}(x)-\kappa_i^+-\frac{2\pi
q(x)}{\kappa|\ln\ep|}\\
\le &\sum_{i=1}^m a^+_{\de,i}\ln\frac R{|x-z^+_i|}/\ln\frac
R{s^+_{\de,i}} -\kappa_i^+
-\frac{2\pi q(x)}{\kappa|\ln\ep|}+o(1)\\
\le &  C\tilde \sigma -\kappa_i^++o(1)<0.
\end{split}
\]

Finally,  if $x\in B_{(s^+_{\de,i})^{\tilde\sigma}}(z^+_i)\setminus
B_{s^+_{\de,i}(1+T(s^+_{\de,i})^{\tilde \sigma})}(z^+_i)$ for some
$i$, then
\[
\begin{split}
&P_{\de,Z}^+(x)-P_{\de,Z}^-(x)-\kappa_i^+-\frac{2\pi
q(x)}{\kappa|\ln\ep|}
= W_{\de,z^+_i,a^+_{\de,i}}(x)-a^+_{\de,i}+O\left(\frac{(s^+_{\de,i})^{\tilde\sigma}}{\ln\frac{R}{s^+_{\de,i}}}\right)\\
=& a^+_{\de,i}\frac{\ln\frac
R{|x-z_i^+|}}{\ln\frac{R}{s^+_{\de,i}}}-a^+_{\de,i}+O\left(\frac{(s^+_{\de,i})^{\tilde\sigma}}{\ln\frac{R}{s^+_{\de,i}}}\right)\\
\le &
-a^+_{\de,i}\frac{\ln(1+T(s^+_{\de,i})^{\tilde\sigma})}{\ln\frac
{R}{s^+_{\de,i}}}+O\left(\frac{(s^+_{\de,i})^{\tilde\sigma}}{\ln\frac{R}{s^+_{\de,i}}}\right)
 <0,
\end{split}
\]
if $T>0$ is large. Note that by the choice of $\tilde\sigma$,
$B_{s^+_{\de,i}(1+(s^+_{\de,i})^\sigma)}(z_i^+)\supset
B_{s^+_{\de,i}(1+T(s^+_{\de,i})^{\tilde\sigma})}(z^+_i)$ for small
$\de$. We therefore derive our conclusion.
\end{proof}

\begin{proposition}\label{ap1}

We have

\[
\begin{split}
I\left(P_{\de,Z}^+-P_{\de,Z}^-\right)=&\frac{C\de^2
 }{\ln\frac{R}{\ep}
 }
 +\frac{\pi(p-1)\de^2}{4(\ln\frac{R}{\ep})^2}\left(\sum_{i=1}^m(\kappa_i^+)^2+\sum_{j=1}^n(\kappa_j^-)^2\right)
 +\sum_{i=1}^m\frac{4\pi^2\de^2\kappa_i^+q(z_i^+)}{|\ln\ep||\ln\frac{R}{\ep}|}\\
 &-\sum_{j=1}^n\frac{4\pi^2\de^2\kappa_j^-q(z_j^-)}{|\ln\ep||\ln\frac{R}{\ep}|}+
 \sum_{i=1}^m \frac{\pi\de^2(\kappa_i^+)^2 g(z_i^+,z^+_i)}{(\ln\frac{R}\ep)^2}+\sum_{j=1}^n \frac{\pi\de^2(\kappa_j^-)^2 g(z_j^-,z^-_j)}{(\ln\frac{R}\ep)^2}
 \\
 &-\sum_{k\ne i}^m \frac{\pi\de^2\kappa_i^+\kappa_k^+ \bar
 G(z_k^+,z^+_i)}{{(\ln\frac{R}\ep)^2}}-\sum_{l\ne j}^n \frac{\pi\de^2\kappa_j^-\kappa_l^- \bar
 G(z_l^-,z^-_j)}{{(\ln\frac{R}\ep)^2}}\\
 &+\sum_{i=1}^m\sum_{j=1}^n\frac{2\pi\de^2\kappa_i^+\kappa_j^-\bar G(z_i^+,z_j^-)}{(\ln\frac{R}{\ep})^2}
 + O\left(\frac{\de^2\ln|\ln\ep|}{|\ln\ep|^3} \right).
 \end{split}
\]

where $C$ is a positive constant.
\end{proposition}

\begin{proof}

Taking advantage of  \eqref{2.3}, we have

\[
\begin{split}
&\de^2\int_\Om \big|D (P_{\de,Z}^+-P_{\de,Z}^-)\big|^2= \sum_{k=1}^m
\sum_{i=1}^m \int_\Om \bigl(W_{\de,z^+_k,a^+_{\de,k}
}-a^+_{\de,k}\bigr)_+^{p}
P_{\de,Z,i}^+\\
&+ \sum_{l=1}^n \sum_{j=1}^n \int_\Om \bigl(W_{\de,z^-_l,a^-_{\de,l}
}-a^-_{\de,l}\bigr)_+^{p} P_{\de,Z,j}^--2 \sum_{j=1}^n \sum_{i=1}^m
\int_\Om \bigl(W_{\de,z^+_i,a^+_{\de,i} }-a^+_{\de,i}\bigr)_+^{p}
P_{\de,Z,j}^-.
\end{split}
\]

First, we estimate

\[
\begin{split}
&\int_{B_{s^+_{\de,i}} (z^+_i)}
\left(W_{\de,z^+_i,a^+_{\de,i}}-a^+_{\de,i}\right)_+^{p}\left(
W_{\de,z^+_i,a^+_{\de,i}}-\frac{a^+_{\de,i}}{\ln\frac
R{s^+_{\de,i}}} g(x,z^+_i)
\right)\\
=&\int_{B_{s^+_{\de,i}} (z^+_i)}
\bigl(W_{\de,z^+_i,a^+_{\de,i}}-a^+_{\de,i}\bigr)^{p+1}+a^+_{\de,i}
\int_{B_{s^+_{\de,i}} (z^+_i)}
\bigl(W_{\de,z^+_i,a^+_{\de,i}}-a^+_{\de,i}\bigr)^{p}\\
&-\frac{a^+_{\de,i}}{\ln\frac R{s^+_{\de,i}}}\int_{B_{s^+_{\de,i}}
(z_i)} \bigl(W_{\de,z^+_i,a^+_{\de,i}}-a^+_{\de,i}\bigr)^{p}
g(x,z_i^+)\\
=&\Bigl(\frac{\de}{s^+_{\de,i}}\Bigr)^{\frac{2(p+1)}{p-1}}(s^+_{\de,i})^2\int_{B_1(0)}\phi^{p+1}
+a^+_{\de,i}\Bigl(\frac{\de}{s^+_{\de,i}}\Bigr)^{\frac{2p}{p-1}}(s^+_{\de,i})^2\int_{B_1(0)}\phi^{p}\\
&-\frac{a^+_{\de,i}}{\ln\frac
R{s^+_{\de,i}}}\Bigl(\frac{\de}{s^+_{\de,i}}\Bigr)^{\frac{2p}{p-1}}g(z^+_i,z^+_i)(s^+_{\de,i})^2\int_{B_1(0)}\phi^{p}+O\left(
\frac{(s^+_{\de,i})^3}{|\ln\ep|^{p+1}}\right)\\
=&\frac{\pi(p+1)}{2}\frac{\de^2(a^+_{\de,i})^2}{(\ln\frac{R}{s^+_{\de,i}})^2}+\frac{2\pi
\de^2(a^+_{\de,i})^2}{\ln\frac{R}{s^+_{\de,i}}}-\frac{2\pi\de^2
(a^+_{\de,i})^2}{(\ln\frac{R}{s^+_{\de,i}})^2}g(z^+_i,z^+_i)+O\left(
\frac{\ep^3}{|\ln\ep|^{p+1}}\right).
\end{split}
\]

Next, for  $k\ne i$,

\[
\begin{split}
&\int_{B_{s^+_{\de,k}} (z^+_k)}
\bigl(W_{\de,z_k^+,a^+_{\de,k}}-a_{\de,k}^+\bigr)_+^{p} \left(
W_{\de,z_i^+,a^+_{\de,i}}-\frac{a^+_{\de,i}}{\ln\frac
R{s^+_{\de,i}}} g(x,z^+_i)
\right)\\
=&
\Bigl(\frac{\de}{s_{\de,k}^+}\Bigr)^{\frac{2p}{p-1}}\frac{a^+_{\de,i}}{\ln\frac{R}{s^+_{\de,i}}}
\int_{B_{s^+_{\de,k}}(z^+_k)}\phi^{p}\Bigl(\frac{|x-z_k^+|}{s_{\de,k}^+}\Bigr)\bar
G(x,z_i^+)
\\
=&\Bigl(\frac{\de}{s^+_{\de,k}}\Bigr)^{\frac{2p}{p-1}}\frac{a^+_{\de,i}(s^+_{\de,k})^2}{\ln\frac{R}{s^+_{\de,i}}}
\bar G(z_k^+,z^+_i)
\int_{B_{1}(0)}\phi^{p}+O\left(\frac{(s^+_{\de,k})^3}{|\ln\ep|^{p+1}}\right)\\
=&\frac{2\pi\de^2
a_{\de,i}^+a_{\de,k}^+}{|\ln\frac{R}{s^+_{\de,i}}||\ln\frac{R}{s^+_{\de,k}}|}\bar{G}(z_i^+,z_k^+)+O\left(
\frac{\ep^3}{|\ln\ep|^{p+1}}\right).
\end{split}
\]
Moreover, we have
\[
\begin{split}
&\int_{B_{s^+_{\de,i}} (z^+_i)}
\bigl(W_{\de,z_i^+,a^+_{\de,i}}-a_{\de,i}^+\bigr)_+^{p} \left(
W_{\de,z_j^-,a^-_{\de,j}}-\frac{a^-_{\de,j}}{\ln\frac
R{s^-_{\de,j}}} g(x,z^-_j)
\right)\\
=&
\Bigl(\frac{\de}{s_{\de,i}^+}\Bigr)^{\frac{2p}{p-1}}\frac{a^-_{\de,j}}{\ln\frac{R}{s^-_{\de,j}}}
\int_{B_{s^+_{\de,i}}(z^+_i)}\phi^{p}\Bigl(\frac{|x-z_i^+|}{s_{\de,i}^+}\Bigr)\bar
G(x,z_j^-)
\\
=&\Bigl(\frac{\de}{s^+_{\de,i}}\Bigr)^{\frac{2p}{p-1}}\frac{a^-_{\de,j}(s^+_{\de,i})^2}{\ln\frac{R}{s^-_{\de,j}}}
\bar G(z_j^-,z^+_i)
\int_{B_{1}(0)}\phi^{p}+O\left(\frac{(s^+_{\de,i})^3}{|\ln\ep|^{p+1}}\right)\\
=&\frac{2\pi\de^2
a_{\de,i}^+a_{\de,j}^-}{|\ln\frac{R}{s^+_{\de,i}}||\ln\frac{R}{s^-_{\de,j}}|}\bar{G}(z_i^+,z_j^-)+O\left(
\frac{\ep^3}{|\ln\ep|^{p+1}}\right).
\end{split}
\]

By Lemma~\ref{al1} and \eqref{2.9},

\[
\begin{split}
&\sum_{k=1}^m\int_\Om \chi_{\Om_k^+}\left( P_{\de,Z}^+-P_{\de,Z}^-
-\kappa_k^+-\frac{2\pi q(x)}{|\ln\ep|}\right)_+^{p+1}\\
=&
\sum_{k=1}^m \int_{B_{Ls^+_{\de,k}}(z_k^+)}\left(
P_{\de,Z}^+-P_{\de,Z}^- -\kappa_k^+-\frac{2\pi
q(x)}{|\ln\ep|}\right)_+^{p+1}
\\
=& \sum_{k=1}^m \int_{B_{Ls^+_{\de,k}}(z_k^+)}
\left(W_{\de,z_k^+,a_{\de,k}^+} -a^+_{\de,k} +O\bigg(\frac
{s_{\de,k}^+}{|\ln\ep|}
\bigg)\right)_+^{p+1}\\
=&\sum_{k=1}^m\left(\frac{\de}{s_{\de,k}^+}\right)^{\frac{2(p+1)}{p-1}}\int_{B_{s^+_{\de,k}}(z^+_k)}\phi^{p+1}\Bigl(\frac{|x-z_k^+|}{s^+_{\de,k}}\Bigr)
+O\left(\frac{(s_{\de,k}^+)^3}{|\ln\ep|^{p+1}}\right)\\
=&\sum_{k=1}^m\left(\frac{\de}{s^+_{\de,k}}\right)^{\frac{2(p+1)}{p-1}}(s_{\de,k}^+)^2\int_{B_{1}(0)}\phi^{p+1}
+O\left(\frac{(s^+_{\de,k})^3}{|\ln\ep|^{p+1}}\right)\\
=&\sum_{k=1}^m\frac{\pi(p+1)}{2}\frac{\de^2(a^+_{\de,k})^2}{(\ln\frac{R}{s^+_{\de,k}})^2}+O\left(
\frac{\ep^3}{|\ln\ep|^{p+1}}\right).
\end{split}
\]
Other terms can be estimated as above.  So, we have proved

\[
\begin{split}
I\left( P_{\de,Z}^+-P_{\de,Z}^-\right)=&\sum_{i=1}^m\left[\frac{\pi
(p+1)}{4} \frac{ \de^2
 (a^+_{\de,i})^2}{|\ln \frac{R}{s^+_{\de,i}}|^2
 } +\frac{\pi\de^2 (a^+_{\de,i})^2}{|\ln \frac{R}{s^+_{\de,i}}|}-\frac{\pi
 g(z_i^+,z_i^+)\de^2(a^+_{\de,i})^2}{|\ln\frac{R}{s^+_{\de,i}}|^2}\right]
 \\
&+\sum_{j=1}^n\left[\frac{\pi (p+1)}{4} \frac{ \de^2
 (a^-_{\de,j})^2}{|\ln \frac{R}{s^-_{\de,j}}|^2
 } +\frac{\pi\de^2 (a^-_{\de,j})^2}{|\ln \frac{R}{s^-_{\de,j}}|}-\frac{\pi
 g(z_j^-,z_j^-)\de^2(a^-_{\de,j})^2}{|\ln\frac{R}{s^-_{\de,j}}|^2}\right]\\
 &+\sum_{k\ne i}^m \frac{\pi\bar G(z_k^+,z^+_i)\de^2 a^+_{\de,i}
 a^+_{\de,k}}{|\ln\frac{R}{s^+_{\de,i}}||\ln \frac{R}{s^+_{\de,k}}|}+\sum_{l\ne j}^n \frac{\pi\bar G(z_l^-,z^-_j)\de^2 a^-_{\de,l}
 a^-_{\de,j}}{|\ln\frac{R}{s^-_{\de,l}}||\ln \frac{R}{s^-_{\de,j}}|}\\
 &-\sum_{i=1}^m\sum_{j=1}^n \frac{2\pi\bar G(z_i^+,z^-_j)\de^2 a^+_{\de,i}
 a^-_{\de,j}}{|\ln\frac{R}{s^+_{\de,i}}||\ln \frac{R}{s^-_{\de,j}}|}
 -\frac{\pi\de^2}{2}\left(\sum_{i=1}^m \frac{(a_{\de,i}^+)^2}{|\ln \frac{R}{s_{\de,i}^+}|^2}\right)\\
 &-\frac{\pi\de^2}{2}\left(\sum_{j=1}^n \frac{(a_{\de,j}^-)^2}{|\ln \frac{R}{s^-_{\de,j}}|^2}\right)
+ O\left(\frac{\ep^3}{|\ln\ep|^{p+1}} \right).
 \end{split}
\]

Thus, the result follows from Remark~\ref{remark2.2}.

\end{proof}

\begin{proposition}\label{ap2}

We have
\[
\begin{split}
\frac{\partial }{\partial z^+_{i,h}}& I\left(
P^+_{\de,Z}-P^-_{\de,Z}\right)
=\frac{4\pi^2\de^2\kappa_i^+}{|\ln\ep||\ln\frac{R}{\ep}|}\frac{\partial
q(z_i^+)}{\partial z^+_{i,h}}
+\frac{2\pi\de^2(\kappa_i^+)^2}{(\ln\frac{R}{\ep})^2} \frac{\partial
g(z_i^+,z_i^+)}{\partial z^+_{i,h}}\\
 &-\sum_{k\ne i}^m \frac{2\pi\de^2\kappa_i^+\kappa_k^+}{(\ln\frac{R}{\ep})^2}\frac{\partial \bar G(z_k^+,z_i^+)}{\partial z^+_{i,h}}
+\sum_{l=1}^n\frac{2\pi\de^2\kappa_i^+\kappa_l^-}{(\ln\frac{R}{\ep})^2}\frac{\partial\bar
G(z_i^+,z_l^-)}{\partial z^+_{i,h}}
+O\left(\frac{\de^2\ln|\ln\ep|}{|\ln\ep|^3} \right),
 \end{split}
\]
\[
\begin{split}
\frac{\partial }{\partial z^-_{j,\bar h}}& I\left(
P^+_{\de,Z}-P^-_{\de,Z}\right)  =
-\frac{4\pi^2\de^2\kappa_j^-}{|\ln\ep||\ln\frac{R}{\ep}|}\frac{\partial
q(z_j^-)}{\partial z^-_{j,\bar h}}
+\frac{2\pi\de^2(\kappa_j^-)^2}{(\ln\frac{R}{\ep})^2} \frac{\partial
g(z_j^-,z_j^-)}{\partial z^-_{j,\bar h}}\\
& -\sum_{l\ne j}^n
\frac{2\pi\de^2\kappa_j^-\kappa_l^-}{(\ln\frac{R}{\ep})^2}\frac{\partial
\bar G(z_l^-,z_j^-)}{\partial z^-_{j,\bar h}}
+\sum_{k=1}^m\frac{2\pi\de^2\kappa_j^-\kappa_k^+}{(\ln\frac{R}{\ep})^2}\frac{\partial\bar
G(z_j^-,z_k^+)}{\partial z^-_{j,\bar h}}
+O\left(\frac{\de^2\ln|\ln\ep|}{|\ln\ep|^3} \right).
 \end{split}
\]

\end{proposition}

\begin{proof}

Direct computation yields that

\[
\begin{split}
\frac{\partial }{\partial z^+_{i,h}} &I\left(
P_{\de,Z}^+-P^-_{\de,Z}\right)\\
=&  \sum_{k=1}^m
\int_{B_{Ls^+_{\de,k}}(z_k^+)}
 \left[ \left(W_{\de,z_k^+,a_{\de,k}^+}-a_{\de,k}^+\right)_+^{p}-
 \left(
P_{\de,Z}^+-P_{\de,Z}^--\kappa_k^+ -\frac{2\pi
q(x)}{|\ln\ep|}\right)_+^{p}\right]\frac{\partial
P_{\de,Z}^+}{\partial z^+_{i,h}}\\
&+ \sum_{l=1}^n \int_{B_{Ls^-_{\de,l}}(z^-_l)}
 \left[ \left(W_{\de,z_l^-,a_{\de,l}^-}-a^-_{\de,l}\right)_+^{p}-
 \left(
P_{\de,Z}^--P_{\de,Z}^+-\kappa_l^- +\frac{2\pi
q(x)}{|\ln\ep|}\right)_+^{p}\right]\frac{\partial
P_{\de,Z}^-}{\partial z^+_{i,h}}\\
& -\sum_{k=1}^m \int_{B_{Ls^+_{\de,k}}(z_k^+)}
 \left[ \left(W_{\de,z_k^+,a_{\de,k}^+}-a_{\de,k}^+\right)_+^{p}-
 \left(
P_{\de,Z}^+-P_{\de,Z}^--\kappa_k^+ -\frac{2\pi
q(x)}{|\ln\ep|}\right)_+^{p}\right]\frac{\partial
P_{\de,Z}^-}{\partial z^+_{i,h}}\\
&-\sum_{l=1}^n \int_{B_{Ls^-_{\de,l}}(z^-_l)}
 \left[ \left(W_{\de,z_l^-,a_{\de,l}^-}-a^-_{\de,l}\right)_+^{p}-
 \left(
P_{\de,Z}^--P_{\de,Z}^+-\kappa_l^+ +\frac{2\pi
q(x)}{|\ln\ep|}\right)_+^{p}\right]\frac{\partial
P_{\de,Z}^+}{\partial z^+_{i,h}}.
 \end{split}
\]

For $k\ne i$, from \eqref{2.9}, we have

\[
\begin{split}
&\int_{B_{Ls_{\de,k}^+}(z_k^+)}
 \left[ \bigl(W_{\de,z_k^+,a_{\de,k}^+}-a_{\de,k}^+\bigr)_+^{p}-
 \left(
P_{\de,Z}^+-P_{\de,Z}^--\kappa_k^+-\frac{2\pi q(x)}{|\ln\ep|} \right)_+^{p}\right]\frac{\partial P^+_{\de,Z,i}}{\partial z^+_{i,h}}\\
=
&\int_{B_{Ls^+_{\de,k}}(z^+_k)}\left[\left(W_{\de,z_k^+,a^+_{\de,k}}-
a_{\de,k}^+\right)^{p-1}\frac{s_{\de,k}^+}{|\ln\ep|}\right] \frac{C}{\ln\frac{R}{s^+_{\de,i}}}\\
 =& O\left(\frac{\ep^3}{|\ln\ep|^{p+1}}
\right).
\end{split}
\]
Using \eqref{2.9}, Lemma~\ref{al1} and Remark~\ref{remark2.2}, we
find that

\[
\begin{split}
&\int_{B_{Ls^+_{\de,i}}(z_i^+)}
 \left[ \left(W_{\de,z^+_i,a^+_{\de,i}}-a^+_{\de,i}\right)_+^{p}-
 \left(
P_{\de,Z}^+-P_{\de,Z}^--\kappa_i^+-\frac{2\pi q(x)}{|\ln\ep|} \right)_+^{p}\right]\frac{\partial P_{\de,Z,i}^+}{\partial z^+_{i,h}}\\
=&\int_{B_{s^+_{\de,i}(1+(s^+_{\de,i})^{\sigma})}(z_i)}
  \left[ \bigl(W_{\de,z^+_i,a^+_{\de,i}}-a^+_{\de,i}\bigr)_+^{p}-
  \left(P_{\de,Z}^+-P_{\de,Z}^--\kappa_i^+-\frac{2\pi q(x)}{|\ln\ep|} \right)_+^{p}\right]\frac{\partial P^+_{\de,Z,i}}{\partial z^+_{i,h}}\\
=
&p\int_{B_{s^+_{\de,i}}(z^+_i)}\bigl(W_{\de,z_i^+,a^+_{\de,i}}-a^+_{\de,i}\bigr)_+^{p-1}\bigg[
 \frac{2\pi}{|\ln\ep|}\bigl\langle D q(z^+_i), x-z^+_i\bigr\rangle+\frac{a^+_{\de,i}}{\ln\frac{R}{s^+_{\de,i}}}\bigl\langle D g(z^+_i,z^+_i), x-z^+_i\bigr\rangle \\
\quad&- \sum_{k\ne i}^m
\frac{a_{\de,k}^+}{\ln\frac{R}{s^+_{\de,k}}}\bigl\langle D \bar
G(z_i^+,z_k^+),x-z_i^+\bigr\rangle +\sum_{l=1}^n
\frac{a_{\de,l}^-}{\ln\frac{R}{s^-_{\de,l}}}\bigl\langle D \bar
G(z_i^+,z_l^-),x-z_i^+\bigr\rangle \bigg]\frac{\partial
P^+_{\de,Z,i}}{\partial z^+_{i,h}}\\
&\,+O\Bigl(\frac{\ep^{2+\sigma}}{|\ln\ep|^{p+1}}\Bigr)\\
=&-\frac{p\de^2a^+_{\de,i}}{|\phi'(1)||\ln\frac{R}{s^+_{\de,i}}|}\bigg(\frac{2\pi}{|\ln\ep|}\frac{\partial
q(z_i^+)}{\partial
z^+_{i,h}}+\frac{a^+_{\de,i}}{\ln\frac{R}{s^+_{\de,i}}}\frac{\partial
g(z_i^+,z_i^+)}{\partial z^+_{i,h}}-\sum_{k\neq
i}^m\frac{a^+_{\de,k}}{\ln\frac{R}{s^+_{\de,k}}}\frac{\partial \bar
G(z_i^+,z_k^+)}{\partial z^+_{i,h}}\\
&+\sum_{l=1}^n\frac{a^-_{\de,l}}{\ln\frac{R}{s^-_{\de,l}}}\frac{\partial
\bar G(z_i^+,z_l^-)}{\partial z^+_{i,h}}\bigg)
\int_{B_1(0)}\phi^{p-1}(|x|)\phi^\prime(|x|)\frac{x_h^2}{|x|}+0\Bigl(\frac{\ep^{2+\sigma}}{|\ln\ep|^{p+1}}\Bigr)\\
 =&\frac{4\pi^2 \de^2a^+_{\de,i}}{|\ln\ep||\ln\frac{R}{s^+_{\de,i}}|}\frac{\partial q(z_i^+)}{\partial z^+_{i,h}}
 +\frac{2\pi \de^2(a^+_{\de,i})^2}{(\ln\frac{R}{s^+_{\de,i}})^2}\frac{\partial g(z_i^+,z^+_i)}{\partial z^+_{i,h}}
 -\sum_{k\neq i}^m\frac{2\pi \de^2a^+_{\de,i}a^+_{\de,k}}{|\ln\frac{R}{s^+_{\de,k}}||\ln\frac{R}{s^+_{\de,i}}|}\frac{\partial \bar G(z^+_i,z_k^+)}{\partial
 z^+_{i,h}}\\
& +\sum_{l=1}^n\frac{2\pi
\de^2a^+_{\de,i}a^-_{\de,l}}{|\ln\frac{R}{s^-_{\de,l}}||\ln\frac{R}{s^+_{\de,i}}|}\frac{\partial
\bar G(z^+_i,z_l^-)}{\partial
 z^+_{i,h}}+O\Bigl(\frac{\ep^{2+\sigma}}{|\ln\ep|^{p+1}}\Bigr),
\end{split}
\]
since

\[
\int_{B_1(0)}\phi^{p-1}(|x|)\phi^\prime(|x|)\frac{x_h^2}{|x|}=
-\frac{2\pi}{p}|\phi^\prime(1)|.
\]

Other terms can be estimated as above.  Thus, the result follows.

\end{proof}

\end{document}